\documentclass{article}

\usepackage{amssymb}
\usepackage{amsmath}
\usepackage{amsthm}
\usepackage{graphicx} 
\usepackage{epstopdf}
\graphicspath{{figures/}}
\usepackage{algorithm}
\usepackage{algorithmic}
\usepackage{multicol}
\usepackage{multirow}
\usepackage{booktabs} 
\usepackage{appendix}
\usepackage{xcolor}
\usepackage{lineno,hyperref}
\usepackage[a4paper]{geometry}
\usepackage{threeparttable}
\newtheorem{theorem}{Theorem}[section]
\newtheorem{proposition}{Proposition}[section]
\newtheorem{lemma}{Lemma}[section]
\newtheorem{definition}{Definition}[section]

\newtheorem{remark}{Remark}
\allowdisplaybreaks[4]
\numberwithin{equation}{section}

\begin{document}
	
	\title{Theory and Fast Learned Solver for $\ell^1$-TV Regularization}
	\author{Xinling Liu\thanks{School of Mathematics $\&$ Statistics, Southwest University, Chongqing {\rm 400715}, China and Key Laboratory of Optimization Theory and Applications at China West Normal University of Sichuan Province, School of Mathematics and Information, China West Normal University, Nanchong {\rm 637009}, China (\texttt{fsliuxl@163.com})} \and Jianjun Wang\thanks{School of Mathematics $\&$ Statistics, Southwest University, Chongqing {\rm 400715}, China (\texttt{wjj@swu.edu.cn, Corresponding author})} \and Bangti Jin\thanks{Department of Mathematics, The Chinese University of Hong Kong, Shatin, New Territories, Hong Kong, China (\texttt{bangti.jin@gmail.com, b.jin@cuhk.edu.hk})}}
	
	
	\date{\today}
	\maketitle
	\begin{abstract}
		The  $\ell^1$ and total variation (TV) penalties have been used successfully in many areas, and the combination of the $\ell^1$ and TV penalties can lead to further improved performance.
		In this work, we investigate the mathematical theory and numerical algorithms for the $\ell^1$-TV model in the context of signal recovery: we derive the sample complexity of the $\ell^1$-TV model for recovering signals with sparsity and gradient sparsity. Also we propose a novel algorithm (PGM-ISTA) for the regularized $\ell^1$-TV problem, and establish its global convergence and parameter selection criteria. Furthermore, we construct a fast learned solver (LPGM-ISTA) by unrolling PGM-ISTA.
		The results for the experiment on ECG signals show the superior performance of LPGM-ISTA in terms of recovery accuracy and computational efficiency.
	\end{abstract}
	
	\section{Introduction}\label{l1tv-intro}
	
	In the last two decades, compressed sensing (CS) has achieved remarkable successes since the pioneering work of Cand\`es, Donoho, Romberg, and Tao \cite{candes2006robust,donoho2006compressed}.
	So far, it has found applications in many fields, e.g., magnetic resonance imaging \cite{lustig2007sparse}, single-pixel imaging \cite{duarte2008signal} and nonlinear parameter identification \cite{JinMaass:2012}.
	The classical setup is to reconstruct an unknown sparse signal $\boldsymbol{x}^\ast\in\mathbb{R}^n$ from linear measurements of the form
	\begin{equation}\label{l1tv-1}
		\boldsymbol{y} = \boldsymbol{A}\boldsymbol{x}^\ast + \boldsymbol{e},
	\end{equation}
	where $\boldsymbol{A}\in\mathbb{R}^{m\times n}(m\ll n)$ is a known sensing matrix and $\boldsymbol{e}\in\mathbb{R}^m$ denotes noise.
	Problem \eqref{l1tv-1} is ill-posed, and the key to address the ill-posedness is to suitably exploit prior structure of real-world signals so as to constrain the solution space. Typical examples include signal sparsity and gradient sparsity, which may appear simultaneously in many real-world datasets.
	This study aims to build a mathematical theory and numerical algorithms of one method for recovering such signals.
	
	First we briefly review single sparsity or gradient sparsity, respectively.
	When the signal is sparse, many mathematical models based on \eqref{l1tv-1} have been proposed, e.g., $\ell^1$ minimization \cite{donoho2006for,wang2020group}, $\ell^p$ ($0\leq p<1$) minimization \cite{wen2015stable,hou2020one,JiaoJinLu:2015} and general nonconvex minimization \cite{wang2019a,HuangJiao:2021}.
	One of the most popular models is the following $\ell^1$ minimization:
	\begin{flalign}\label{l1tv-2}		\qquad&\textup{($\ell^1$):}\qquad\qquad\qquad\qquad\min_{\boldsymbol{x}\in\mathbb{R}^n}\|\boldsymbol{x}\|_1,\quad \textup{s.t.}\quad \|\boldsymbol{A}\boldsymbol{x}-\boldsymbol{y}\|_2\leq\epsilon,&
	\end{flalign}
	where $\epsilon\geq0$ denotes noise level (i.e. $\|\boldsymbol{e}\|_2\leq\epsilon$).
	If $\boldsymbol{x}^\ast$ is $s$-sparse (i.e. $\|\boldsymbol{x}^\ast\|_0\leq s$, where $\|\boldsymbol{x}^\ast\|_0$ denotes the number of nonzero entries of $\boldsymbol{x}^\ast$) and $\boldsymbol{A}$ satisfies certain properties (e.g., null sparse property and restricted isometry property), then $\boldsymbol{x}^\ast$ can be effectively recovered from a few measurements by, e.g., orthogonal matching pursuit (OMP) \cite{mo2012a,wen2017a} and iterative shrinkage-threshold algorithm (ISTA) \cite{daubechies2004an} or its fast version (FISTA) \cite{beck2009a}.
	Signals with few details are often locally constant with a few jumps, which can be well represented as sparse gradients \cite{hou2024tensor,liu2025guaranteed}. A signal $\boldsymbol{x}^\ast$ with sparse gradients takes the form of a sparse $\boldsymbol{D}\boldsymbol{x}^\ast$ (with $\boldsymbol{D}=\operatorname{circ}([-1,1,0,\cdots,0])\in \mathbb{R}^{(n-1)\times n}$ denoting a row circulant difference matrix).
	One popular approach to recover such $\boldsymbol{x}^\ast$ is to solve
	\begin{flalign}\label{l1tv-3}
		\qquad&\textup{(TV):}\qquad\qquad\qquad\qquad\min_{\boldsymbol{x}\in\mathbb{R}^n}\|\boldsymbol{x}\|_{\rm TV},\quad \textup{s.t.}\quad \|\boldsymbol{A}\boldsymbol{x}-\boldsymbol{y}\|_2\leq\epsilon,&
	\end{flalign}
	where $\|\boldsymbol{x}\|_{\rm TV}=\|\boldsymbol{D}\boldsymbol{x}\|_1$ denotes the total variation (TV) seminorm of $\boldsymbol{x}$.
	Cai et al. \cite{cai2015guarantees} proved that the TV method can recover a gradient $s$-sparse signal $\boldsymbol{x}^\ast$ (i.e., $\|\boldsymbol{D}\boldsymbol{x}^\ast\|_0\leq s$) if the sampling number $m$ obeys $m=\mathcal{O}(\sqrt{sn}\ln n)$ and the entries of  $\boldsymbol{A}$ follow the standard Gaussian distribution.
	However, the bound is not optimal, and several works have improved the result  \cite{daei2018sample,daei2019on,genzel2021l1,genzel2022compressed}.
	The TV problem \eqref{l1tv-3} can be effectively solved using, e.g., proximal gradient descent algorithm and projected subgradient algorithm \cite{beck2017first}.
	
\begin{figure}[t]
		\centering
		{\includegraphics[width=\textwidth]{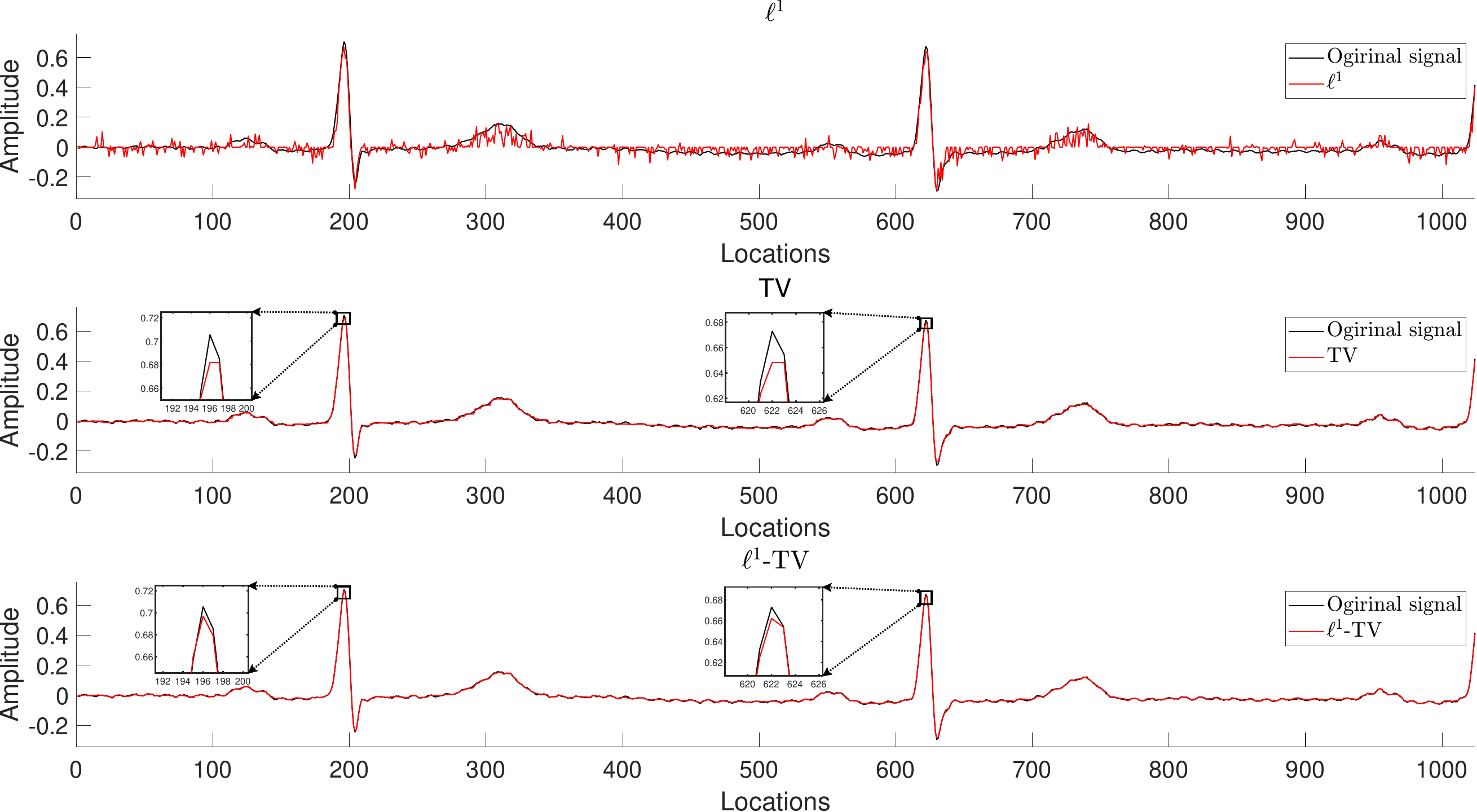}}
		\caption{The comparison of $\ell^1$, TV and $\ell^1$-TV under Gaussian sampling with a sampling ratio of 0.5.}
		\label{l1tv-intro-show}
	\end{figure}

	The use of the $\ell^1$ and TV penalties promotes respectively sparsity and gradient sparsity in the signals.
	However, one single penalty cannot well characterize signals with both properties, e.g., electrocardiogram (ECG) signals \cite{garudari2010artifacts} and time-varying signals  \cite{angelosante2009compressed}.
	Tibshirani et al. \cite{tibshirani2005sparsity} proposed to combine the $\ell^1$ and TV penalties to promote sparsity in both signals and their successive differences, which has achieved great success in many applications \cite{tang2017fused,cui2021fused,mohammadi2019a,wang2015fused,xin2014efficient}, e.g., regression, feature selection
	and denoising. However, it has not been studied in the context of CS so far. This motivated the theoretical study of the following problem:
	\begin{flalign}\label{l1tv-model}
		\qquad&(\ell^1\text{-}{\rm TV}):\qquad\qquad\qquad\qquad \min_{\boldsymbol{x}\in\mathbb{R}^n} \lambda_1\|\boldsymbol{x}\|_1 + \lambda_2\|\boldsymbol{x}\|_{\rm TV}, \quad \textup{s.t.} \quad \|\boldsymbol{A}\boldsymbol{x}-\boldsymbol{y}\|_2\leq\epsilon,&
	\end{flalign}
	where regularization parameters $\lambda_1\geq0$ and $\lambda_2\geq0$.
	Obviously, if $\lambda_1=0$ and $\lambda_2>0$, \eqref{l1tv-model} reduces to \eqref{l1tv-2}, and if $\lambda_1>0$ and $\lambda_2=0$, \eqref{l1tv-model} becomes \eqref{l1tv-3}.
	To illustrate the benefit of the $\ell^1$-TV method, we take one ECG signal from the MIT-BIH Arrhythmia Database \cite{moody2001the} as the original signal and apply random Gaussian sampling with a ratio of 0.5 and $\epsilon=0$. The restorations by the $\ell^1$, TV and $\ell^1$-TV methods are shown in Fig. \ref{l1tv-intro-show}, which shows that the restorations by the $\ell^1$ and TV methods are unsatisfactory, whereas $\ell^1$-TV gives a very satisfactory reconstruction.
	Thus, the $\ell^1$-TV model can enforce sparsity in both signals and the differences	between consecutive components of the signals in the CS framework. To theoretically substantiate this empirical observation, we shall establish the required sampling number for the robust recovery of signals regarding given sparsity and gradient sparsity levels.
	
	The second important issue is the algorithmic aspect of the model \eqref{l1tv-model}.
	The model \eqref{l1tv-model} poses big challenges to the numerical solution since it involves both non-separable term $\|\boldsymbol{x}\|_{\rm TV}$ and non-smooth term $\|\boldsymbol{x}\|_1$, as well as	quadratic inequality constraint.
	Li et al. \cite{li2018on} suggested a level-set method for solving \eqref{l1tv-model}, which however can be computationally very challenging due to the complicated geometry of the quadratic inequality constraint, and the potentially unknown noise level.
	In this work, we consider the regularized model
	\begin{equation}\label{l1tv-un}
		\qquad\min_{\boldsymbol{x}\in\mathbb{R}^n} \tfrac{1}{2}\|\boldsymbol{y}-\boldsymbol{A}\boldsymbol{x}\|_2^2 + \lambda_1\|\boldsymbol{x}\|_1 + \lambda_2\|\boldsymbol{D}\boldsymbol{x}\|_1,
	\end{equation}
	by slightly abusing the notation $\lambda_1$ and $\lambda_2$.
	Note that such a strategy is common in the community, and there are several works on the equivalence between the constrained problem and its regularized counterpart \cite{li2022selecting,anzengruber2009morozov,wen2011parameter}.
	Transforming problem \eqref{l1tv-model} to \eqref{l1tv-un} avoids the need to estimate the noise level $\epsilon$ and facilitates implementing existing algorithms, e.g., linearized alternating direction method of multipliers (LADMM) \cite{li2014linearized} and smooth-FISTA (S-FISTA) \cite{beck2012smoothing}.
	
	To circumvent these issues, we develop an algorithm unrolling approach for problem \eqref{l1tv-un} such that all algorithmic hyper-parameters in the corresponding algorithm (not including the regularization parameters $\lambda_1$ and $\lambda_2$) are learnable. Algorithm unrolling is a powerful technique: it inherits the interpretability of iterative algorithms and high expressivity of deep neural networks \cite{monga2021algorithm}.
	Since its first introduction by Gregor et al. \cite{gregor2010learning}, the idea has found many successful applications, e.g., signal and imaging processing, including compressed sensing \cite{yang2018admm}, sparse coding \cite{gregor2010learning} and tomographic reconstruction \cite{BarbanoKereta:2022}; see the the review \cite{monga2021algorithm} for details.
	In algorithm unrolling, each iteration typically corresponds to one layer of a neural network, and the employed operator must be adapted to the activation function used by the neural network.
	For example, in ISTA (for the $\ell^1$ minimization), the soft-threshold operator is the activation function in learned ISTA (LISTA) \cite{gregor2010learning}.
	Similarly, the 1D TV minimization can also be effectively unfolded into a recurrent neural network (i.e., LPGD-Taut) \cite{cherkaoui2020learning}.
	These works motivates developing an algorithm unrolling approach for solving problem \eqref{l1tv-un}.
	However, due to the presence of two penalties, directly applying these important findings to \eqref{l1tv-un} is infeasible.
	We shall develop a novel fast learned solver for \eqref{l1tv-un} by unrolling a new algorithm with guaranteed global convergence and rigorous parameter selection criteria.
	
	The main contributions of the work are as follows.
	\begin{itemize}
		\item We present novel theoretical guarantees on the robust recovery via \eqref{l1tv-model}.
		The bound on the sampling number reflects the joint influence of sparsity and gradient sparsity of the signal, where higher levels of sparsity or gradient sparsity result in a larger sampling number.
		\item We propose a novel algorithm, i.e., proximal gradient mapping ISTA (PGM-ISTA) for solving \eqref{l1tv-un}, and prove its global convergence and parameter selection.
		Furthermore, we develop an efficient learned solver LPGM-ISTA for the model \eqref{l1tv-un} based on algorithm unrolling.
		\item We conduct extensive numerical experiments, including ECG signals.
		The results show that LPGM-ISTA already perform very well with only 2 layers/iterations, and outperforms traditional algorithms in terms of computing time and accuracy.
	\end{itemize}
	
	Lastly, we situate the present work in existing works.
	First, in the analysis, to obtain the robust recovery guarantee of the $\ell^1$-TV method, we select a vector from the subdifferential of the objective function in \eqref{l1tv-model} that yields an upper bound on its statistical dimension, and then derive a bound on the sampling number for successful signal recovery.
	Several studies \cite{candes2011compressed,kabanava2015analysis,kabanava2015robust,genzel2021l1} have addressed closely related topics.
	These studies can provide recovery guarantee for the model \eqref{l1tv-model}, but they focus on analysis sparsity, and thus do not deliver the sampling number for the sparsity and gradient sparsity levels of the concerned signal.
	Second, PGM-ISTA is based on proximal operator and gradient mapping operator \cite{beck2017first}, which directly solves problem \eqref{l1tv-un}, rather than its equivalent or approximate alternatives.
	Thus, PGM-ISTA may have a lower computational complexity per iteration, compared to LADMM \cite{li2014linearized}, which requires the update of auxiliary variables.
	Moreover, the S-FISTA \cite{beck2012smoothing} is an approximate approach for solving problem \eqref{l1tv-un}, and its performance depends heavily on the smoothing parameter.
	Traditional optimization algorithms for problem \eqref{l1tv-un} are often time-consuming since they require many iterations (dozens to hundreds).
	The unrolled network LPGM-ISTA provides an efficient way to tackle the computational challenge.
	
	The rest of this paper is structured as follows.
	In Section \ref{l1tv-section2}, we recall useful notation, and some technical lemmas, and
	in Section \ref{l1tv-section3}, present the robust guarantee for the $\ell^1$-TV method \eqref{l1tv-model}.
	In Section \ref{l1tv-section4}, we propose a learned solver by unrolling the PGM-ISTA algorithm, and
	in Section \ref{l1tv-simulation}, provide numerical illustrations. All the proofs are given
	in Section \ref{l1tv-proof}. We conclude in Section \ref{l1tv-section6} with additional discussions.
	
	\section{Preliminaries}\label{l1tv-section2}
	
	First we recall useful notation. Throughout, scalars are denoted by lowercase case letters (e.g. $x$), vectors by bold lowercase letters (e.g. $\boldsymbol{x}$), and matrices by bold uppercase letters (e.g. $\boldsymbol{A}$).
	The $i$-th element of a vector $\boldsymbol{x}$ is denoted by $x_i$ or $[\boldsymbol{x}]_i$.
	$\boldsymbol{0}$ denotes a zero vector, and	$\boldsymbol{I}_n$ the identity matrix of size $n\times n$.
	We denote  sets by hollow capitalized letters (e.g. $\mathbb{S}$), and denote by $|\mathbb{S}|$ its cardinality.
	$\mathbb{N}_+$ refers to the set of positive integers, and $[n]=\{1,2,\cdots,n\}$ and $\mathbb{S}^c$ the complement of $\mathbb{S}$.
	For a vector $\boldsymbol{x}\in\mathbb{R}^n$, $\boldsymbol{x}_\mathbb{S}$ is the vector in $\mathbb{R}^n$ that coincides with $\boldsymbol{x}$ on the entries in $\mathbb{S}$ and zero otherwise.
	For a vector $\boldsymbol{x}\in\mathbb{R}^n$ and $p>0$, we denote $\|\boldsymbol{x}\|_p=\big(\sum_{i=1}^{n}|x_i|^p)\big)^{1/p}$ as the $\ell^p$ norm of $\boldsymbol{x}$, and we also denote $\|\boldsymbol{x}\|_0 = |\{i:x_i \neq 0\}|$.
	For a matrix $\boldsymbol{A}$, we denote its spectral norm by $\|\boldsymbol{A}\|_2$.
	For any $\boldsymbol{x}\in\mathbb{R}^n$ and a set $\mathbb{S}\subset\mathbb{R}^n$, we denote the distance of $\boldsymbol{x}$ and $\mathbb{S}$ by $\textup{dist}(\boldsymbol{x},\mathbb{S})=\min_{\boldsymbol{z}\in\mathbb{S}}\|\boldsymbol{x}-\boldsymbol{z}\|_2$.
	The operator $[a]_+=\max\{a,0\}$ gives the positive part of the scalar $a$.
	For $\boldsymbol{x},\boldsymbol{y}\in\mathbb{R}^n$, $\boldsymbol{x}\odot\boldsymbol{y}$ denotes the Hadamard / componentwise product.
	We denote the support of $\boldsymbol{x}\in\mathbb{R}^n$ by $\mathbb{S}_R(\boldsymbol{x}):{=}\{i\in[n]:x_i\neq 0\}$ ($\mathbb{S}_R$ for short) and the gradient support by $\mathbb{S}_G(\boldsymbol{x}):{=}\{i\in[n-1]:x_{i+1}-x_i\neq 0\}$ ($\mathbb{S}_G$ for short).
	We denote by $\mathbb{S}^\ast$ the set of solutions to problem \eqref{l1tv-un}.
	
	First we recall the definitions of descent cone, subdifferential, and statistical dimension.
	Recall that a proper convex function takes at least one finite value but never the value $-\infty$.
	We denote the set of  extended real numbers by $\overline{\mathbb{R}}=\mathbb{R}\cup\{\pm\infty\}$.
	
	\begin{definition}[{\cite{tropp2015convex}}]
		Let $g:\mathbb{R}^n\rightarrow\overline{\mathbb{R}}$ be a proper convex function.
		The descent cone $\mathcal{D}(g,\boldsymbol{x})$ of the function $g$ at a point $\boldsymbol{x}\in\mathbb{R}^n$ is defined by
		\begin{equation*}
			\mathcal{D}(g,\boldsymbol{x}) = \mathop{\cup}\limits_{t>0}\left\{\boldsymbol{u}\in\mathbb{R}^n:g(\boldsymbol{x}+t\boldsymbol{u})\leq g(\boldsymbol{x})\right\}.
		\end{equation*}
	\end{definition}
	
	\begin{definition}[{\cite{rockafellar2015convex}}]
		The subdifferential of a convex function $f:\mathbb{R}^n\rightarrow\mathbb{R}$ at a point $\boldsymbol{x}\in \mathbb{R}^n$ is given by
		$\partial f(\boldsymbol{x}) = \{\boldsymbol{z}\in\mathbb{R}^n:f(\boldsymbol{y})\geq f(\boldsymbol{x}) + \langle \boldsymbol{z},\boldsymbol{y}-\boldsymbol{x}\rangle,\ \forall\boldsymbol{y}\in\mathbb{R}^n\}$.
	\end{definition}
	
	\begin{definition}[{\cite{rockafellar2015convex}}]
		For a given non-empty set $\mathbb{K}\subset\mathbb{R}^n$, the cone obtained by $\mathbb{K}$ is defined by $ \operatorname{cone}({\mathbb{K}})=\{\lambda\boldsymbol{x}\in\mathbb{R}^n,\boldsymbol{x}\in \mathbb{K},\lambda\geq0\}$.
	\end{definition}
	
	\begin{definition}[{\cite{amelunxen2014living}}]
		Let $\mathbb{K}\subset\mathbb{R}^n$ be a convex closed cone and $\boldsymbol{g}\sim \mathcal{N}(\boldsymbol{0},\boldsymbol{I}_n)$ be a standard Gaussian vector in $\mathbb{R}^n$. The statistical dimension of $\mathbb{K}$ is defined by
		\begin{equation}\label{l1tv-6}
			\delta(\mathbb{K}) = \operatorname{E}[ \textup{dist}^2(\boldsymbol{g},\mathbb{K}^\circ)]=\min_{\boldsymbol{x}\in\mathbb{K}^\circ}\|\boldsymbol{g}-\boldsymbol{x}\|_2,
		\end{equation}
		where $\operatorname{E}[\cdot]$ denotes the operator for expectation and $\mathbb{K}^\circ=\{\boldsymbol{z}\in\mathbb{R}^n:\langle\boldsymbol{z},\boldsymbol{v}\rangle\leq0\ \forall \boldsymbol{v}\in \mathbb{K}\}$ denotes the set of outward normals of $\mathbb{K}$, i.e., the polar cone of $\mathbb{K}$.
	\end{definition}
	
	The descent cone is connected with the subdifferential by \cite{rockafellar2015convex}
	\begin{equation}\label{l1tv-7}
		[\mathcal{D}(g,\boldsymbol{x})]^\circ = \textup{cone} (\partial g(\boldsymbol{x})) = \bigcup\limits_{t\geq0}\{t\boldsymbol{z}: \boldsymbol{z}\in\partial g(\boldsymbol{x})\}.
	\end{equation}
	Combining \eqref{l1tv-6} and \eqref{l1tv-7} gives
	\begin{align}\label{l1tv-8}
		\delta(\mathcal{D}(g,\boldsymbol{x})) &= \operatorname{E}[\textup{dist}^2(\boldsymbol{g},\mathcal{D}(g,\boldsymbol{x})^\circ)]= \operatorname{E}[\textup{dist}^2(\boldsymbol{g},\textup{cone}(\partial g(\boldsymbol{x})))]\nonumber\\
		&=\operatorname{E}\Big[\inf_{\substack{t\geq0\\ \boldsymbol{z}\in\partial g(\boldsymbol{x})}}\ \|\boldsymbol{g}-t\boldsymbol{z}\|_2^2\Big].
	\end{align}
	
	Next, we give the definition of conic Gaussian width and its relation with statistical dimension.
	\begin{definition}[{\cite{tropp2015convex}}]
		Let $\mathbb{S}^{n-1}$ be the $(n-1)$-dimensional unit sphere, and  $\mathbb{K}\subset\mathbb{R}^n$ be a cone, not necessarily convex. The conic Gaussian width $w(\mathbb{K})$ is defined by
		\begin{equation}\label{l1tv-25}
			w(\mathbb{K}) = \operatorname{E}\Big[\sup_{\boldsymbol{u} \in \mathbb{K}\cap \mathbb{S}^{n-1}}\langle\boldsymbol{g},\boldsymbol{u}\rangle\Big].
		\end{equation}
	\end{definition}
	
	The conic Gaussian width is useful for estimating the sampling number $m$.
	\begin{lemma}[{\cite[Proposition 10.2]{amelunxen2014living}}]
		\label{l1tv-lemma1}
		Let $\mathbb{K}$ be a convex cone. Then the statistical dimension $\delta(\mathbb{K})$ and the conic Gaussian width $w(\mathbb{K})$ satisfy
		\begin{equation}
			w^2(\mathbb{K})\leq\delta(\mathbb{K})\leq w^2(\mathbb{K}) + 1.
		\end{equation}
	\end{lemma}
	
	Combining \eqref{l1tv-8} with Lemma \ref{l1tv-lemma1} gives
	\begin{equation}\label{l1tv-k16}
		w^2(\mathcal{D}(g,\boldsymbol{x})) \leq \operatorname{E}\Big[\inf_{\substack{t\geq0\\ \boldsymbol{z}\in\partial g(\boldsymbol{x})}}\ \|\boldsymbol{g}-t\boldsymbol{z}\|_2^2\Big].
	\end{equation}
	
	\begin{definition}[{\cite[Definition 5.1]{beck2017first}}]\label{l1tv-lsmooth}
		Fix $l\geq0$. A function $f:\mathbb{R}^n\rightarrow (-\infty,\infty]$ is said to be  $l$-smooth over a set $\mathbb{K}\subset\mathbb{R}^n$ if it is differentiable over $\mathbb{K}$ and satisfies
		$\|\nabla f(\boldsymbol{x})-\nabla f (\boldsymbol{y})\|_2\leq l\|\boldsymbol{x}-\boldsymbol{y}\|_2$ for all $\boldsymbol{x}, \boldsymbol{y} \in \mathbb{K}$,
		and	the constant $l$ is called the smoothness parameter.
	\end{definition}
	
	\begin{definition}[{\cite[Definition 6.1]{beck2017first}}]
		\label{l1tv-proxmap}
		For the function $f:\mathbb{K}\rightarrow(-\infty,\infty]$, the proximal mapping of $f$ at $\boldsymbol{x}\in\mathbb{K}$ is the operator given by
		\begin{align*}
			\operatorname{prox}_f(\boldsymbol{x})=\mathop{\arg\min}_{\boldsymbol{u}\in\mathbb{K}} f(\boldsymbol{u}) + \tfrac{1}{2}\|\boldsymbol{u}-\boldsymbol{x}\|_2^2.
		\end{align*}
	\end{definition}
	
	For $f(\boldsymbol{x})=\lambda_1\|\boldsymbol{x}\|_1$, the proximal mapping is given by the soft threshold operator:
	\begin{align}\label{l1tv-102}
		\mathcal{S}_{\lambda_1}(\boldsymbol{x}):{=}\operatorname{prox}_f(\boldsymbol{x})=\operatorname{sign}(\boldsymbol{x})\odot\max\{0,|\boldsymbol{x}|-\lambda_1\},
	\end{align}
	For $f(\boldsymbol{x})=\lambda_2\|\boldsymbol{x}\|_{\textup{TV}}$, we denote $\mathcal{T}_{\lambda_2}(\boldsymbol{x}):{=}\operatorname{prox}_f(\boldsymbol{x})$. It can  be computed via dynamic programming, e.g., taut-string algorithm \cite{condat2013a} at an $\mathcal{O}(n)$ complexity.
	Thus, $\mathcal{P}_{\lambda_1}^{\lambda_2}(\boldsymbol{x})=\operatorname{prox}_g(\boldsymbol{x})$ can be computed  efficiently.
	\begin{lemma}[{\cite[Theorem 1]{liu2010an}}]
		\label{l1tv-lemma4}
		Let $\boldsymbol{x}\in\mathbb{R}^n$ and $g(\boldsymbol{x})=\lambda_1\|\boldsymbol{x}\|_1+\lambda_2\|\boldsymbol{x}\|_{\textup{TV}}$. Then there holds
		\begin{align*}
			\mathcal{P}_{\lambda_1}^{\lambda_2}(\boldsymbol{x}):{=}\operatorname{prox}_g(\boldsymbol{x}) =\mathcal{S}_{\lambda_1}(\mathcal{T}_{\lambda_2}(\boldsymbol{x})).
		\end{align*}
	\end{lemma}

	\section{Recovery guarantee for $\ell^1$-TV regularization}\label{l1tv-section3}
	
	Now we provide robust recovery guarantee of the model \eqref{l1tv-model}.
	We first give Theorem \ref{l1tv-theorem3.1} for estimating the statistical dimension $\delta(\mathcal{D}(g,\boldsymbol{x}))$ of $\mathcal{D}(g,\boldsymbol{x})$.
	It is crucial for proving Theorem \ref{l1tv-theorem3.2} on robust recovery.
	Below we denote the regular sparsity by $s_r(\boldsymbol{x})=|\mathbb{S}_R(\boldsymbol{x})|$ ($s_r$ for short), and the gradient sparsity by $s_g(\boldsymbol{x})=|\mathbb{S}_G(\boldsymbol{x})|$ ($s_g$ for short). 
	All the proofs are given in Section \ref{l1tv-proof1}.
	\begin{theorem}\label{l1tv-theorem3.1}
		Let $n\in\mathbb{N}_+$, ${\boldsymbol{x}}\in\mathbb{R}^n$, and $s_r=|\mathbb{S}_R|$,  $s_g=|\mathbb{S}_G|$.
		Then we have
		\begin{align}
			\delta(\mathcal{D}(g,{\boldsymbol{x}}))\leq \Phi(s_r,s_g),
		\end{align}
		with
		\begin{align}\label{l1tv-phi}
			\Phi(s_r,s_g)=n-\frac{6}{\pi}\frac{\left[\lambda_1(n-s_r)+\sqrt{2}\lambda_2(n-1-s_g)\right]^2}{3n\lambda_1^2+4\left(2n+s_g-4\right)\lambda_2^2+12\lambda_1\lambda_2\min\{s_r,s_g\}}.
		\end{align}
	\end{theorem}
	
	\begin{remark}\label{l1tv-remark11}
		Theorem \ref{l1tv-theorem3.1} gives an upper bound on $\delta(\mathcal{D}(g,{\boldsymbol{x}}))$, which is crucial for the recovery guarantee of \eqref{l1tv-model}. By Lemma \ref{l1tv-lemma1} and the estimate \eqref{l1tv-k16}, we can bound  the conic Gaussian width of $\mathcal{D}(g,\boldsymbol{x})$ by
		\begin{align}\label{l1tv-83}
			w^2(\mathcal{D}(g,\boldsymbol{x}))\leq\delta(\mathcal{D}(g,{\boldsymbol{x}}))\leq\Phi(s_r,s_g).
		\end{align}
		If $\lambda_1=0$, the $\ell^1$-TV method reduces to the TV method, and \eqref{l1tv-83} implies
		\begin{align*}
			w^2(\mathcal{D}(\|\cdot\|_{\textup{TV}},\boldsymbol{x}))\leq\Phi(s_r,s_g)=n-\frac{3}{\pi}\frac{(n-s_g-1)^2}{2n+s_g-4}\overset{\triangle}{=}\Phi_{\rm TV},
		\end{align*}
		which is identical with that in \cite[Theorem 1]{daei2018sample} for the TV method. If $\lambda_2=0$, the $\ell^1$-TV method reduces to the $\ell^1$ method, and \eqref{l1tv-83} implies
		\begin{align}\label{l1tv-84}
			w^2(\mathcal{D}(\|\cdot\|_1,\boldsymbol{x}))\leq\Phi(s_r,s_g)=n-\frac{2}{\pi}\frac{(n-s_r)^2}{n}.
		\end{align}
		This bound is also meaningful, in view of the order sharp one, i.e. $\Phi_{\ell^1}=2s_r\log(n/s_r)+2s_r$ \cite{tropp2015convex}.
		More precisely, \eqref{l1tv-84} always represents a sharper upper bound than that in \cite{tropp2015convex} for a signal with length $n=10^5$ when the sparsity level $s_r<5000$.
		Since conic Gaussian width is a typical tool for deriving lower bound for successful recovery, \eqref{l1tv-84} may provide a tighter sampling bound for the $\ell^1$ method.
	\end{remark}
	
	\begin{remark}
		We further examine how the parameters $n$, $\lambda_1$, $\lambda_2$ influence the variation of $\Phi(s_r,s_g)$ with respect to the sparsity levels $s_r$ and $s_g$.
		Since $\Phi(\cdot,\cdot)$ only depends on the ratio of $\lambda_1$ and $\lambda_2$, we may fix $\lambda_2=1$.
		For multiple combinations of $n$ and $\lambda_1$, since $s_g\leq2s_r$ is always satisfied, we plot the function values of $\Phi(s_r,s_g)$ in the region $\{(s_r,s_g)|s_g\leq2s_r\}$ in Fig. \ref{l1tv-srsg}.
		The following observations can be drawn from the plots: (1) $\Phi(s_r,s_g)$ exhibits consistent patterns across different choices of $n$ and $\lambda_1$;
		(2) $\Phi(s_r,s_g)$ is not sensitive to $s_r$, and remains relatively stable for fixed $s_g$, indicating that $\Phi(s_r,s_g)$ can dynamically reflect the influence of sparsity and gradient sparsity levels for estimating conic Gaussian width to some extent;
		(3) For fixed $s_g$, $\Phi(s_r,s_g)$ shows a monotonic increase pattern with respect to $s_r$, which may still make sense for certain small $s_r$ according to Remark \ref{l1tv-remark11}.
		\begin{figure*}[htp]
			\renewcommand{\arraystretch}{0.5}
			\setlength\tabcolsep{0.5pt}
			\centering
			\begin{tabular}{ccccccc}
				\centering
				\includegraphics[width=40mm, height = 30mm]{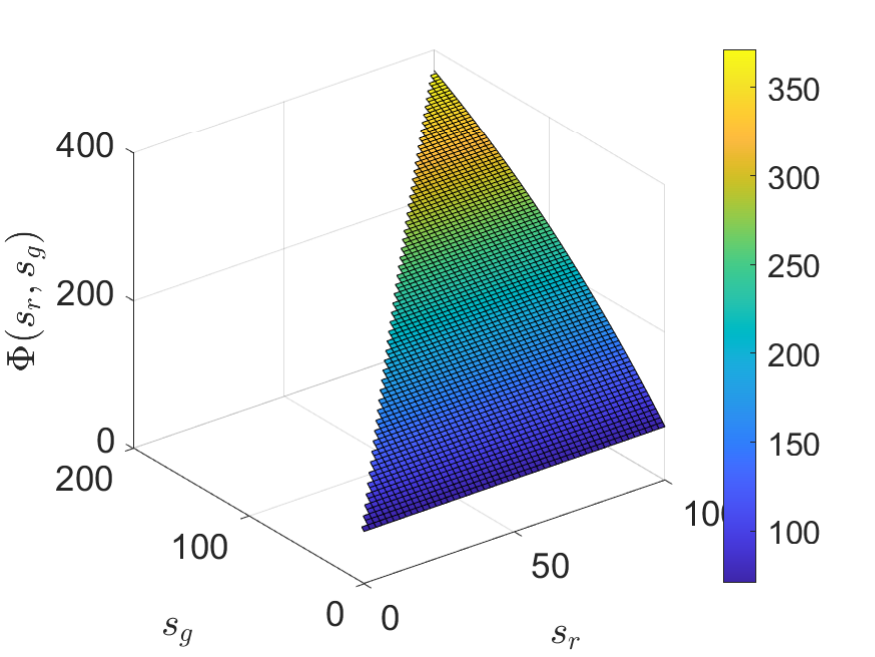}&
				\includegraphics[width=40mm, height = 30mm]{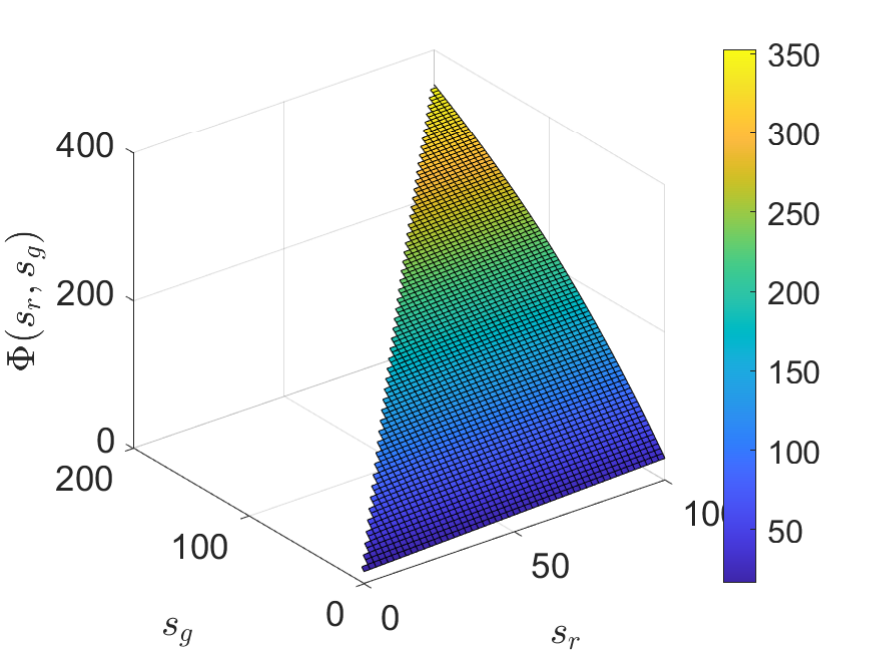}&
				\includegraphics[width=40mm, height = 30mm]{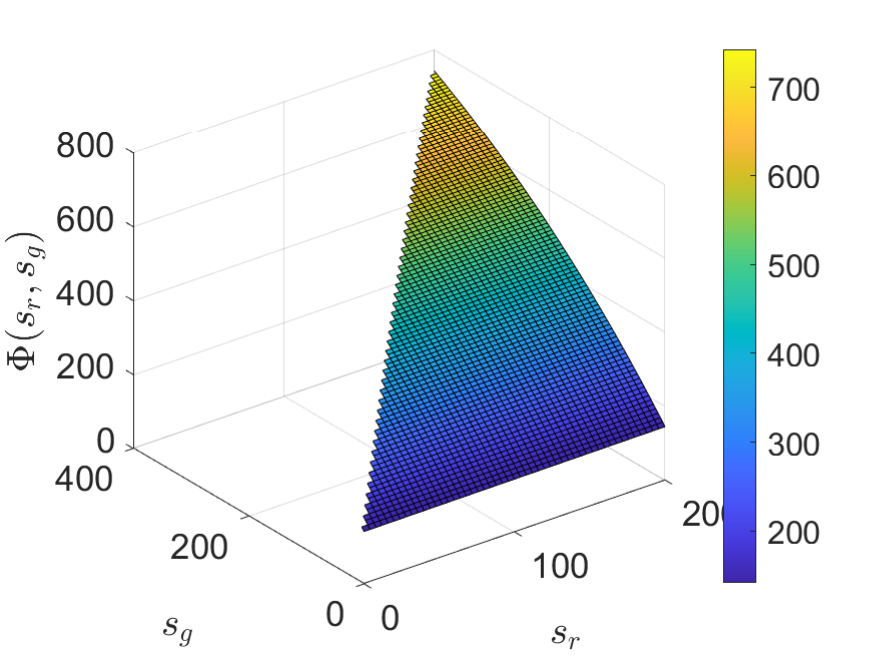}&
				\includegraphics[width=40mm, height = 30mm]{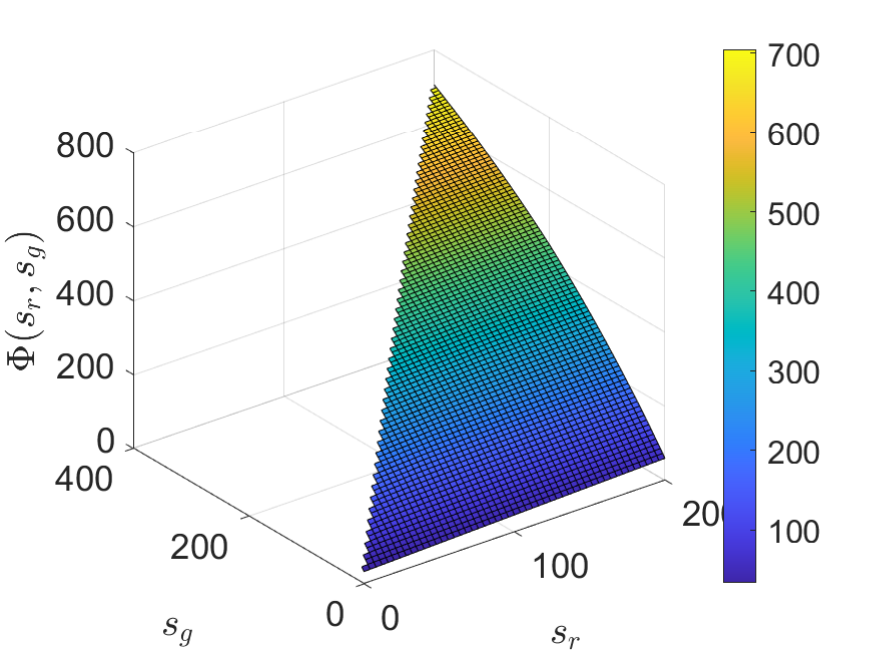}\\
				\scriptsize $n=500$, $\lambda_1=0.01$&
				\scriptsize $n=500$, $\lambda_1=0.1$&
				\scriptsize $n=1000$, $\lambda_1=0.01$&
				\scriptsize $n=1000$, $\lambda_1=0.1$
			\end{tabular}
			\caption{For fixed $\lambda_2 = 1$, plots of $\Phi(s_r,s_g)$ ($s_g\leq2s_r$) under multiple settings of $n$ and $\lambda_1$.}
			\label{l1tv-srsg}
		\end{figure*}
	\end{remark}
	
	\begin{remark}\label{l1tv-rem3}
		Let $n_0$ be an upper bound on the conic Gaussian width of the $\ell^1$ or TV methods. We examine the case $\Phi(s_r,s_g)\leq n_0$ more closely.
		Direct computation shows
		\begin{align}\label{l1tv-101}
			2(n-s_r)^2\Big(\frac{\lambda_1}{\lambda_2}\Big)^2+b\frac{\lambda_1}{\lambda_2}+c\geq0,
		\end{align}
		with $b=4\sqrt{2}(n-s_r)(n-1-s_g)+(3n+12\min\{s_r,s_g\})n_0$ and $c=4(n-1-s_g)^2+4(2n+s_g-4)n_0$.
		\eqref{l1tv-101} is a quadratic inequality in $\lambda_1/\lambda_2$, whose solution in $\lambda_1/\lambda_2>0$ is always nonempty. That is, a proper choice of $n$, $s_r$, $s_g$ and $\lambda_2/\lambda_1$ can always lead to a solution of \eqref{l1tv-101}. To shed insights, we compare $\Phi(s_r,s_g)$ with $\Phi_{\ell^1}$ for $\ell^1$ method \eqref{l1tv-2} and $\Phi_{\rm TV}$ for TV method \eqref{l1tv-3} in Table \ref{l1tv-tab1}.
		The table indicates that when $\lambda_1/\lambda_2=1$, $\Phi(s_r,s_g)$ can be much smaller than $\Phi_{\ell^1}$ and $\Phi_{\rm TV}$, and in the relatively extreme case $\lambda_1/\lambda_2=0.1$ (i.e., gradient sparsity is dominating), the presence of sparsity can generally still improve $\Phi(s_r,s_g)$ over $\Phi_{\rm TV}$.
		
		\begin{table}[hbt!]
			\setlength{\tabcolsep}{1.8mm}
			\centering
			\fontsize{7}{8}\selectfont
			\begin{threeparttable}
				\caption{Comparisons of the upper bounds of conic Gaussian width for $\ell^1$, TV and $\ell^1$-TV methods with $n=1000$.}
				\label{l1tv-tab1}
				\begin{tabular}{ccccccccccccccc}
					\toprule[1pt]
					\multirow{2}{*}{Methods}&
					\multicolumn{2}{c}{$s_r=50$}&
					\multicolumn{2}{c}{$s_r=100$}&
					\multicolumn{2}{c}{$s_r=150$}\cr
					\cline{2-7}
					&$s_g=25$&$s_g=50$&$s_g=50$&$s_g=100$&$s_g=75$&$s_g=150$\cr
					\hline
					\hline
					&\multicolumn{6}{c}{$\lambda_1/\lambda_2=1$}\cr
					\hline
					\hline
					$\Phi_{\ell^1}$\cite{tropp2015convex}&400&400&661&661&870&870\cr
					$\Phi_{\rm TV}$\cite{daei2018sample}&552&580&580&632&670&680\cr
					$\Phi(s_r,s_g)$&92&149&186&285&271&400\cr
					\hline
					\hline
					&\multicolumn{6}{c}{$\lambda_1/\lambda_2=0.1$}\cr
					\hline
					\hline
					$\Phi_{\ell^1}$\cite{tropp2015convex}&400&400&660&661&870&870\cr
					$\Phi_{\rm TV}$\cite{daei2018sample}&552&580&580&632&607&680\cr
					$\Phi(s_r,s_g)$&508&553&556&632&599&697\cr
					\bottomrule[1pt]
					\vspace{-0.8cm}
				\end{tabular}
			\end{threeparttable}
		\end{table}	
	\end{remark}
	
	Using Theorem \ref{l1tv-theorem3.1}, we can give a recovery guarantee of the $\ell^1$-TV method.
	See Section \ref{l1tv-proof2} for the proof.
	\begin{theorem}\label{l1tv-theorem3.2}
		Let $\boldsymbol{x}^\ast\in\mathbb{R}^n$, $\epsilon\geq0$, $\boldsymbol{A}\in\mathbb{R}^{m\times n}$ whose rows are independent random vectors drawn from the standard Gaussian distribution $\mathcal{N}(\boldsymbol{0},\boldsymbol{I}_n)$, and let $\boldsymbol{y} = \boldsymbol{A}\boldsymbol{x}^\ast + \boldsymbol{e}\in\mathbb{R}^m$ be the measurements. Suppose that $\|\boldsymbol{e}\|_2\leq \epsilon$ and $\hat{\boldsymbol{x}}$ is the solution of the $\ell^1\textup{-TV}$ method \eqref{l1tv-model}.
		For certain $t>0$, if
		\begin{equation}\label{l1tv-85}
			m>(\sqrt{\Phi(s_r,s_g)}+t)^2+1,
		\end{equation}
		where $\Phi(s_r,s_g)$ is defined in \eqref{l1tv-phi} with $s_r(\boldsymbol{x})$ and $s_g(\boldsymbol{x})$ replaced by $s_r(\boldsymbol{x}^\ast)$ and $s_g(\boldsymbol{x}^\ast)$, then, with a probability at least $1-e^{-\frac{t^2}{2}}$, there holds
		\begin{equation*}
			\|\boldsymbol{x}^\ast - \hat{\boldsymbol{x}}\|_2\leq \frac{2\epsilon}{\sqrt{m-1}-\Phi(s_r,s_g)-t}.
		\end{equation*}
	\end{theorem}
	
	\begin{remark}
		Theorem \ref{l1tv-theorem3.2} gives robust recovery guarantee of the model \eqref{l1tv-model}. The sampling number for successful recovery depends on the regular sparsity level $s_r$ and gradient sparsity level $s_g$ of the signal $\boldsymbol{x}^*$ to be recovered.
		By Remark \ref{l1tv-remark11}, we can also obtain  the sampling numbers of the $\ell^1$ method and the TV method. If $\lambda_1=0$, the sampling number of the TV method defined by \eqref{l1tv-85} is consistent with that from \cite{daei2018sample}.
		If $\lambda_2=0$, the sampling number for the robust recovery of the $\ell^1$ method defined by \eqref{l1tv-85} is $m>\big(\sqrt{n-\frac{2}{\pi}\frac{(n-s_r)^2}{n}}+t\big)^2+1$, which is sharper than $m>\big(\sqrt{2s_r\log(n/s_r)+2s_r}+t\big)^2+1$ in \cite{tropp2015convex} for certain $n$ and $s_r$.
		Moreover, Remark \ref{l1tv-rem3} indicates that the sample complexity of the $\ell^1$-TV method can be better than both $\ell^1$ \cite{tropp2015convex} and TV \cite{daei2018sample} for suitable regimes for $\lambda_1/\lambda_2$, $s_r$ and $s_g$.
	\end{remark}
	
	\begin{remark}
		Note that we can simplify $\lambda_1\|\boldsymbol{x}\|_1+\lambda_2\|\boldsymbol{D}\boldsymbol{x}\|_1$ as $\|\Psi\boldsymbol{x}\|_1$, where $\Psi=[\lambda_1\boldsymbol{I}_n;\lambda_2\boldsymbol{D}]\in\mathbb{R}^{(2n-1)\times n}$ is a frame with frame bounds $\lambda_1$ and $\sqrt{\lambda_1^2+4\lambda_2^2}$\footnote{For any $\boldsymbol{x}\in\mathbb{R}^n$, the inequalities $\lambda_1\|\boldsymbol{x}\|_2\leq\|\Psi\boldsymbol{x}\|_2\leq\sqrt{\lambda_1^2+4\lambda_2^2}\|\boldsymbol{x}\|_2$ can be easily checked by using the estimate  $\|\boldsymbol{D}\|_2\leq2$ \cite[Lemma 8.11]{zhang2023numerical}.}.
		The studies \cite{candes2011compressed,kabanava2015analysis,kabanava2015robust} have considered the case of $\Psi$ being a general frame.
		Cand\`es et al. \cite{candes2011compressed} studied the case based on D-RIP (RIP adapted to $\Psi$), which is NP-hard to be verified.
		Kabanava et al. \cite{kabanava2015analysis,kabanava2015robust} studied the construction of recovery condition via tangent cones.
		Genzel et al. \cite{genzel2021l1} derived an explicit formula that describes the precise required number of measurements, whose bound depends on the coherence structure of $\Psi$.
		However, all of these studies do not treat the influence of the sparsity and gradient sparsity levels of the signal on the sampling number needed for robust recovery.
	\end{remark}

	\section{Fast learned solver for regularized $\ell^1$-TV model}\label{l1tv-section4}
	
	The model \eqref{l1tv-un} is numerically challenging  due to its nonsmoothness and nonseparability.
	Existing algorithms, e.g., LADMM and S-FISTA, have shortcomings, e.g., dependence on additional hyper-parameters.
	Moreover, the iterative approach often suffers from relatively high cost. To mitigate these issues, we propose a novel algorithm and study its convergence, and then construct a fast learned solver by unrolling the proposed algorithm. Throughout, we also rewrite the model \eqref{l1tv-un} as
	\begin{equation}\label{l1tv-gen}
		\min_{\boldsymbol{x}\in\mathbb{R}^n} F(\boldsymbol{x}) = f(\boldsymbol{x}) + g(\boldsymbol{x}),
	\end{equation}
	with $f(\boldsymbol{x}) = \frac{1}{2}\|\boldsymbol{y}-A\boldsymbol{x}\|_2^2$, $g(\boldsymbol{x})=g_1(\boldsymbol{x}) + g_2(\boldsymbol{x})$ with $g_1(\boldsymbol{x}) = \lambda_1\|\boldsymbol{x}\|_1$ and $g_2(\boldsymbol{x}) =\lambda_2 \|\boldsymbol{x}\|_{\textup{TV}}$.
	
	\subsection{Proximal gradient mapping method}\label{ssec:PGM}
	
	First we construct a proximal gradient mapping method, prove its global convergence and discuss the parameter selection criterion.
	
	\subsubsection{Construction of proximal gradient mapping method}
	
	We first recall the proximal gradient descent method and gradient mapping operator, which are crucial for constructing the proximal gradient mapping method.
	The goal of the proximal gradient descent method \cite{beck2017first} is to find a minimum of the function $f + g$, where $f$ is convex and smooth (with a Lipschitz constant $L_f$), and $g$ is convex, not necessarily smooth.
	It consists of the following iterative step
	\begin{equation}\label{l1tv-15}
		\boldsymbol{x}^{k+1} = \operatorname{prox}_{\frac{1}{L_f}g}\left(\boldsymbol{x}^k - \tfrac{1}{L_f}\nabla f(\boldsymbol{x}^k)\right).
	\end{equation}
	When $g(\boldsymbol{x}) = \|\boldsymbol{x}\|_1$, the iterative step \eqref{l1tv-15} has a closed form solution (see \eqref{l1tv-102}), resulting in ISTA or FISTA.
	However, this family of methods cannot be directly employed to solve \eqref{l1tv-un} with $g(\boldsymbol{x}) = \lambda_1\|\boldsymbol{x}\|_1 + \lambda_2\|\boldsymbol{x}\|_{\textup{TV}}$.
	In fact, since $g(\cdot)$ is a sum of $\ell^1$ composite terms, computing $\operatorname{prox}_g(\cdot)$ is no longer separable and ISTA loses its effectiveness.
	In addition, the gradient mapping operator of $f$ and $g$ is given by
	\begin{equation*}
		\mathcal{G}_{L_f}^{f,g}(\boldsymbol{x}) = L_f\left(\boldsymbol{x} - \operatorname{prox}_{\frac{1}{L_f}g}\left(\boldsymbol{x} - \tfrac{1}{L_f}\nabla f(\boldsymbol{x})\right)\right).
	\end{equation*}
	Then the ISTA update step in \eqref{l1tv-15} can be rewritten as $\boldsymbol{x}^{k+1} = \boldsymbol{x}^k - \frac{1}{L_f}\mathcal{G}_{L_f}^{f,g}(\boldsymbol{x}^k)$, which takes the form of a gradient descent step.
	$\mathcal{G}_{L_f}^{f,g}(\cdot)$ generalizes the gradient of $f(\boldsymbol{x})$:
	\begin{equation*}
		\begin{cases}
			\mathcal{G}_{L_f}^{f,g}(\boldsymbol{x}) = \nabla h(\boldsymbol{x}) = \nabla f(\boldsymbol{x}),& \text{if $g(\boldsymbol{x}) \equiv 0$},\\
			\mathcal{G}_{L_f}^{f,g}(\boldsymbol{x}^\ast) = 0, &\text{if and only if $\boldsymbol{x}^\ast$ is a minimizer of $f+g$}.
		\end{cases}
	\end{equation*}
	See \cite{beck2017first} for more details about the gradient mapping operator.
	
	Now we can develop a proximal gradient mapping (PGM) method.
	For problem \eqref{l1tv-un}, motivated by the proximal gradient descent method and gradient mapping operator, we propose to solve \eqref{l1tv-gen} using a PGM method, which takes an update of the form
	\begin{equation*}
		\boldsymbol{x}^{k+1} = \operatorname{prox}_{tg_2}\left(\boldsymbol{x}^k - t\mathcal{G}_{1/u}^{f,g_1}(\boldsymbol{x}^k)\right),
	\end{equation*}
	where $u>0$ and $t>0$ are constants to be specified. This update involves computing $\operatorname{prox}_{ug_1}(\cdot)=\mathcal{S}_{\lambda_1u}(\boldsymbol{x})$ and $\operatorname{prox}_{tg_2}(\cdot)=\mathcal{T}_{\lambda_2t}(\boldsymbol{x})$ (see \eqref{l1tv-gen} for $g_1(\cdot)$ and $g_2(\cdot)$), both of which can be carried out efficiently.
	Thus, the proposed PGM method with the inner update step being the soft threshold shrinkage and the outer update step being the proximal TV (PGM-ISTA for short) can be concisely written as
	\begin{equation*}
		\boldsymbol{x}^{k+1} = \mathcal{T}_{\lambda_2t}\Big(\big(1-\tfrac{t}{u}\big)\boldsymbol{x}^k+\tfrac{t}{u}\mathcal{S}_{\lambda_1u}\big(\boldsymbol{x}^k - u\boldsymbol{A}^{\top}(\boldsymbol{A}\boldsymbol{x}^k - \boldsymbol{y})\big)\Big).
	\end{equation*}
	The whole procedure of PGM-ISTA is given in Algorithm \ref{alg:pgm-ista}.
	PGM-ISTA involves two adjustable parameters (i.e., $u$ and $t$) due to the two parameters in the $\ell^1$-TV objective.
	
	\begin{algorithm}[hbt!]
		\caption{PGM$-$ISTA}
		\label{alg:pgm-ista}
		\begin{algorithmic}
			\STATE {\textbf{Input}:  $(f,g_1,g_2,\boldsymbol{x}^0)$}, where $f$, $g_1$ and $g_2$ satisfy \eqref{l1tv-gen}, and $\boldsymbol{x}^0\in\mathbb{R}^n$.
			\STATE {\textbf{Initialization}: pick $u>0$, $t>0$}.
			\STATE {\textbf{General Step}: for $k=0,1,2,\cdots$, excute the following steps:}\\
			$(a)$ set $\mathcal{G}_{1/u}^{f,g_1}(\boldsymbol{x}^k)=\frac{1}{u}\left(\boldsymbol{x}^k-\operatorname{prox}_{ug_1}\left(\boldsymbol{x}^k-u\nabla f(\boldsymbol{x}^k)\right)\right)$;\\
			$(b)$ set $\boldsymbol{x}^{k+1}=\operatorname{prox}_{tg_2}\left(\boldsymbol{x}^k-t\mathcal{G}_{1/u}^{f,g_1}(\boldsymbol{x}^k)\right)$.
		\end{algorithmic}
	\end{algorithm}
	
	\subsubsection{Global convergence and parameter selection criterion for PGM-ISTA}\label{l1tv-sec4.1.1}
	
	It is natural to ask whether PGM-ISTA converges to the optimal solution of \eqref{l1tv-un}.
	Now we discuss its global convergence and parameter selection criterion.
	The proof of the next result is given in Section \ref{l1tv-proof3}.
	\begin{theorem}\label{l1tv-thm2}
		Let $\{\boldsymbol{x}^k(u,t)\}_{k\geq0}$ be a sequence generated by Algorithm \ref{alg:pgm-ista}. Then the following statements hold.
		\begin{itemize}
			\item[{\rm(a)}] There exist $u>0$ and $t>0$ such that $\{\boldsymbol{x}^k(u,t)\}_{k\geq0}$ converges to the optimal solution of problem \eqref{l1tv-un};
			\item[{\rm(b)}] Suppose $u\in(0,2/\|\boldsymbol{A}\|_2^2)$, $t\in(0,u]$.
			Then, $\{\boldsymbol{x}^k(u,t)\}_{k\geq0}$ converges to a fixed point $\boldsymbol{x}(u,t)$:
			\begin{equation}
				\boldsymbol{x}(u,t)=\operatorname{prox}_{tg_2}\big(\boldsymbol{x}(u,t)-t\mathcal{G}_{1/u}^{f,g_1}\left(\boldsymbol{x}(u,t)\right)\big).
			\end{equation}
		\end{itemize}
	\end{theorem}
	
	\begin{remark}\label{l1tv-remark1}
		The proof of Theorem \ref{l1tv-thm2} {\rm(a)} indicates that if the parameter pair $(u_0,t_0)$ provides a convergent sequence generated by Algorithm \ref{alg:pgm-ista}, smaller parameter pair $(u,t)$ with $u\leq u_0$ and some $t>0$ also generates a convergent sequence.
	\end{remark}
	
	\begin{remark}\label{l1tv-rem1}
		Theorem \ref{l1tv-thm2}{\rm(b)} indicates that for $u\in(0,2/\|\boldsymbol{A}\|_2^2)$ and $t\in(0,u]$, the sequence generated by Algorithm \ref{alg:pgm-ista} always converges to some fixed point $\boldsymbol{x}^\infty$, that is
		\begin{equation}\label{l1tv-20}
			\boldsymbol{x}^\infty = \operatorname{prox}_{tg_2}\Big(\big(1-\tfrac{t}{u}\big)\boldsymbol{x}^\infty + \tfrac{t}{u}\hat{\boldsymbol{x}}\Big),
		\end{equation}
		with $\hat{\boldsymbol{x}}=\operatorname{prox}_{ug_1}\left(\boldsymbol{x}^\infty - u\nabla f(\boldsymbol{x}^\infty)\right)$.
		Therefore, by \cite[Theorem 6.39]{beck2017first}, \eqref{l1tv-20} holds if and only if
		\begin{equation*}
			\exists\  \boldsymbol{w}_1\in \partial g_1(\hat{\boldsymbol{x}}), \boldsymbol{w}_2\in\partial g_2(\boldsymbol{x}^\infty)\quad \textup{s.t. } \boldsymbol{w}_1 + \boldsymbol{w}_2 +\nabla f(\boldsymbol{x}^\infty) = \boldsymbol{0}.
		\end{equation*}
		By Fermat's optimality condition \cite[Theorem 3.63]{beck2017first}, $\boldsymbol{x}^\ast \in \mathbb{S}^\ast$ is an optimal solution of \eqref{l1tv-gen} if and only if
		\begin{equation*}
			\exists\  \boldsymbol{w}_1\in \partial g_1(\boldsymbol{x}^\ast), \boldsymbol{w}_2\in\partial g_2(\boldsymbol{x}^\ast)\quad  \textup{s.t. \ } \boldsymbol{w}_1 + \boldsymbol{w}_2 +\nabla f(\boldsymbol{x}^\ast) = \boldsymbol{0}.
		\end{equation*}
		These two observations imply that a sufficient condition for the fixed point $\boldsymbol{x}^\infty$ to be the optimal solution of \eqref{l1tv-un} is $\hat{\boldsymbol{x}} = \boldsymbol{x}^\infty$.
		These discussions indicate that $\boldsymbol{x}^\infty$ may not be the optimal solution of \eqref{l1tv-gen}.
	\end{remark}
	
	Now we discuss the parameter selection strategy. Theorem \ref{l1tv-thm2} not only shows that there always {\color{red}exists} a parameter pair $(u,t)$ for Algorithm \ref{alg:pgm-ista} to generate a convergent sequence to the optimal solution of \eqref{l1tv-gen}, but also provides guidelines for choosing the parameter pair $(u,t)$.
	Remark \ref{l1tv-rem1} indicates that although a fixed point of Algorithm \ref{alg:pgm-ista} may not necessarily be an optimal solution of \eqref{l1tv-gen}, it narrows down the feasible range since the optimal solution must be a fixed point.
	Naturally, if we define the set of fixed points of Algorithm \ref{alg:pgm-ista} by
	\begin{equation}\label{l1tv-23}
		\mathbb{F} = \left\{\boldsymbol{x}(u,t)\in\mathbb{R}^n|u\in\left(0,\tfrac{2}{\|\boldsymbol{A}\|_2^2}\right),t\in(0,u]\right\},
	\end{equation}
	then Theorem \ref{l1tv-thm2} implies $\mathbb{S}^\ast\cap\mathbb{F}\neq\emptyset$ (empty set).
	
	The next theorem further investigates the influence of the parameter pair $(u,t)$ on the convergence of Algorithm \ref{alg:pgm-ista}, provided that $g_1(\cdot)$ has a bounded domain $\operatorname{dom}(g_1) \subset\{\boldsymbol{x} : \|\boldsymbol{x}\|_2 \leq r\}$ for sufficiently large $r$. Thus, for any $\boldsymbol{x} \in \operatorname{dom}(g_1)$, the gradient $\nabla f$ of $f$ can be bounded by
	\begin{equation}\label{l1tv-5}
		\|\nabla f(\boldsymbol{x})\|_2 \leq r\|\boldsymbol{A}\|_2^2 + \|\boldsymbol{A}^{\top}\boldsymbol{y}\|_2.
	\end{equation}
	The proof is deferred to Section \ref{l1tv-proof4}.
	\begin{theorem}\label{l1tv-thm4}
		Let $u\in(0,2/\|\boldsymbol{A}\|_2^2)$, $t\in(0,3u/4)$. Then for any $\epsilon > 0$, if $k\in\mathbb{N}$ is sufficiently large such that
		\begin{equation}\label{l1tv-86}
			\|{\boldsymbol{x}}^k - {\boldsymbol{x}}^{k+1}\|_2\leq t\epsilon,
		\end{equation}
		we have
		\begin{equation*}
			F(\boldsymbol{x}^{k+1}) - F(\boldsymbol{x}^\ast) \leq 2r\epsilon +(t\kappa + \beta)u,
		\end{equation*}
		with $\boldsymbol{x}^\ast = \arg\min_{\boldsymbol{x}} F(\boldsymbol{x})$ and
		\begin{align*}
			&\kappa = \left(\|\boldsymbol{A}\|_2^2r+\|\boldsymbol{A}^{\top}\boldsymbol{y}\|_2+\lambda_2\|\boldsymbol{D}\|_2+\lambda_1\right)\|\boldsymbol{A}\|_2^2(\|\boldsymbol{A}\|_2^2r + \|\boldsymbol{A}^{\top}\boldsymbol{y}\|_2+\lambda_1),\\
			&\beta = 2r\|\boldsymbol{A}\|_2^2(\|\boldsymbol{A}\|_2^2r + \|\boldsymbol{A}^{\top}\boldsymbol{y}\|_2+\lambda_1) + (\|\boldsymbol{A}\|_2^2r+\|\boldsymbol{A}^{\top}\boldsymbol{y}\|_2)(\|\boldsymbol{A}\|_2^2r+\lambda_1)+\tfrac{\lambda_1^2}{2}.
		\end{align*}
	\end{theorem}
	
	\begin{remark}\label{l1tv-remark2}
		Theorem \ref{l1tv-thm4} shows that for $u\in(0,2/\|\boldsymbol{A}\|_2^2)$ and $t\in(0,3u/4)$, the error in the function value is controlled by $2r\epsilon +(t\kappa + \beta)u$, and thus a smaller $u$ will lead to smaller errors in the function value for any fixed $\epsilon$.
		This and Theorem \ref{l1tv-thm2} show that we need smaller $u$ and proper $t$ to get a smaller objective function error in practice, i.e., a balance between $u$ and $t$ in order to achieve a better rate of convergence.
	\end{remark}
	
	\subsection{Construction of fast learned solver}
	In this part, we construct a learned solver for the model \eqref{l1tv-un} by unrolling PGM-ISTA. The discussions in Section \ref{ssec:PGM} suggest that the PGM-ISTA scheme can provide an efficient iterative solver for the model \eqref{l1tv-un}.
	By setting the first iteration of PGM-ISTA to $\boldsymbol{x}^0=\boldsymbol{0}$, it becomes
	\begin{equation}\label{l1tv-24}
		\boldsymbol{x}^1=\mathcal{T}_{t\lambda_2}(\tfrac{t}{u}\mathcal{S}_{u\lambda_1}(u\boldsymbol{A}^{\top}\boldsymbol{y})).
	\end{equation}
	Formally this can viewed as a neural network with two layers. Next, we present more evidences to elucidate the point.
	Indeed, by setting proper learnable parameters, the ISTA with the soft threshold operator $\mathcal{S}_\lambda$ can be extended to	a recurrent network, i.e., learned ISTA (LISTA) \cite{gregor2010learning}.
	Thus, the inner part of \eqref{l1tv-24} (i.e., $\mathcal{S}_{u\lambda_1}(\cdot)$) can be regarded as one layer of a neural network similar to LISTA.
	Further, the study on the TV proximal operator \cite{cherkaoui2020learning} shows that the outer part of \eqref{l1tv-24} (i.e. $\mathcal{T}_{t\lambda_2}(\cdot)$) can also be learned. These useful findings are summarized in the following proposition.
	\begin{proposition}[{\cite[Weak Jacobian of prox-TV]{cherkaoui2020learning}}]\label{l1tv-thm5}
		Let $\boldsymbol{x}\in\mathbb{R}^n$ and $\boldsymbol{z}=\mathcal{T}_{\mu}(\boldsymbol{x})$, and  $\mathbb{S}$ be the support of $\widetilde{\boldsymbol{D}}\boldsymbol{z}$ with $\widetilde{\boldsymbol{D}}=[[1,0,\cdots,0];\boldsymbol{D}]$.
		Then, the weak Jacobian $J_{\boldsymbol{x}}$ and $J_\mu$ of $\mathcal{T}_{\mu}(\boldsymbol{x})$ with respect to $\boldsymbol{x}$ and $\mu$ are given by \begin{equation*}
			\left\{\begin{aligned}		 J_{\boldsymbol{x}}(\boldsymbol{x},\mu)&=\boldsymbol{L}_{:,\mathbb{S}} (\boldsymbol{L}_{:,\mathbb{S}}^{\top}\boldsymbol{L}_{:,\mathbb{S}})^{-1}\boldsymbol{L}_{:,\mathbb{S}}^{\top},\\
				J_\mu(\boldsymbol{x},\mu)&=-\boldsymbol{L}_{:,\mathbb{S}}(\boldsymbol{L}_{:,\mathbb{S}}^{\top}\boldsymbol{L}_{:,\mathbb{S}})^{-1}
				\operatorname{sign}(\boldsymbol{D}\boldsymbol{z})_{\mathbb{S}},
			\end{aligned}\right.
			\quad\mbox{with }
			\boldsymbol{L}=\left[
			\begin{array}{cccccc}
				1  & 0 & \cdots & 0\\
				1 & 1 & \cdots & 0\\
				\vdots & \vdots & \ddots & \vdots\\
				1 & 1 & \cdots & 1
			\end{array}
			\right]\in \mathbb{R}^{n\times n}.
		\end{equation*}
	\end{proposition}
	
	Proposition \ref{l1tv-thm5} gives the weak Jacobians of $\mathcal{T}_\mu(\boldsymbol{x})$ respect to $\boldsymbol{x}$ and $\mu$, which can be employed via back-propagation.
	Moreover, the dependency in the inputs is only through $\mathbb{S}$ and $\operatorname{sign}(\boldsymbol{D}\boldsymbol{z})$.
	Hence, computing these weak Jacobians can be done efficiently by simply storing $\operatorname{sign}(\boldsymbol{D}\boldsymbol{z})$ as a mask, as it would be done for a ReLU or the soft-thresholding activation functions, and requiring just $2(n-1)$ bits.
	
\begin{figure}[hbt!]
		\centering
		{\includegraphics[width=\textwidth]{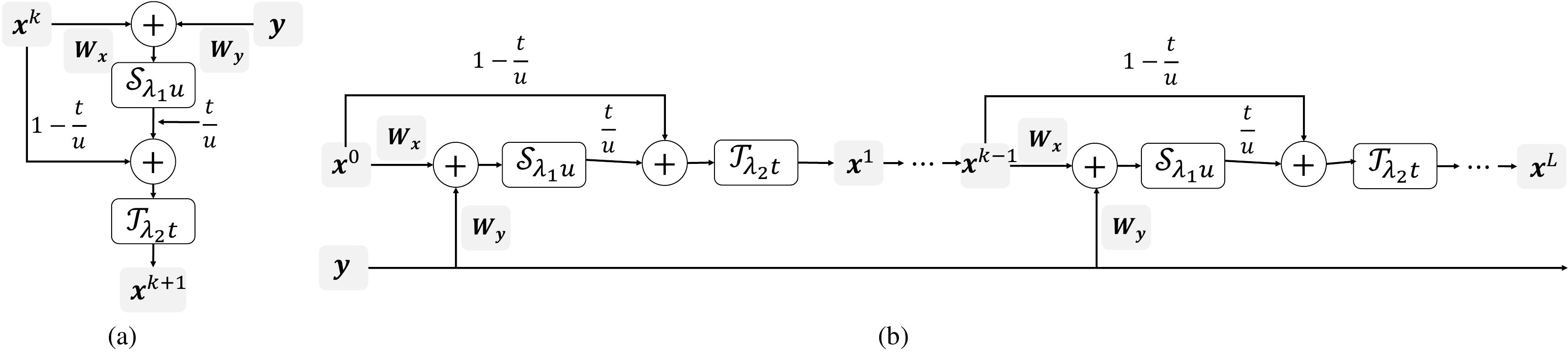}}
		\caption{Neural network representation of PGM-ISTA. (a) PGM-ISTA; (b) LPGM-ISTA obtained by unfolding PGM-ISTA with $L$ layers, and the learnable parameters are $\Theta^L=\{\boldsymbol{W}_{\boldsymbol{x}}, \boldsymbol{W}_{\boldsymbol{y}},u,t\}$.}
		\label{network_unfolding}
	\end{figure}

	Now one may interpret some neural networks as approximate algorithms for solving \eqref{l1tv-un}.
	To this end, let $\boldsymbol{W}_{\boldsymbol{y}} = u\boldsymbol{A}^{\top}$, $\boldsymbol{W}_{\boldsymbol{x}} = \boldsymbol{I}_n-u\boldsymbol{A}^{\top}\boldsymbol{A}$.
	Then, a general iteration of PGM-ISTA takes the form
	\begin{equation}\label{l1tv-iteration}
		\boldsymbol{x}^{k+1} = \mathcal{T}_{\lambda_2t}\left(\left(1-\tfrac{t}{u}\right)\boldsymbol{x}^k+\tfrac{t}{u}\mathcal{S}_{\lambda_1u}\big(\boldsymbol{W}_{\boldsymbol{x}}\boldsymbol{x}^k+\boldsymbol{W}_{\boldsymbol{y}}\boldsymbol{y}\big)\right).
	\end{equation}
	Fig. \ref{network_unfolding} (a) depicts PGM-ISTA for problem \eqref{l1tv-un}, and the theoretical findings in Section \ref{l1tv-sec4.1.1} ensure the convergence of iteration \eqref{l1tv-iteration} with properly chosen $u\in(0,2/\|\boldsymbol{A}\|_2^2)$ and $t\in(0,u]$.
	By unrolling the architecture with $L$ layers, we obtain a network $net(\boldsymbol{y};{\Theta^L})=\boldsymbol{x}^L$ which is termed as learned PGM-ISTA (LPGM-ISTA) with learnable parameters $\Theta^L=\{\boldsymbol{W}_{\boldsymbol{x}}, \boldsymbol{W}_{\boldsymbol{y}},u,t\}$, cf. Fig. \ref{network_unfolding} (b) for a schematic illustration.
	
	By Remark \ref{l1tv-remark2}, a good choice of $u$ and $t$ is crucial in order to achieve fast convergence.
	In practice, the parameter selection scheme from the theory of PGM-ISTA can be used to initialize the  training stage.
	The goal of the learning process is to find the parameter vector $\Theta^L$ that minimizes the loss or the expected loss over a suitable training data distribution,	which is however unavailable.
	The empirical risk minimization framework with a given training set $\{(\boldsymbol{y}_i,{\boldsymbol{x}}_i)\}_{i=1}^{N}$ is a practical choice.
	Thus, the neural network is trained by minimizing the following loss function
	\begin{equation}
		\min_{\Theta^L}\frac{1}{N}\sum_{i=1}^{N}\|\boldsymbol{x}_i-net(\boldsymbol{y}_i;{\Theta^L})\|_2^2.
	\end{equation}
	Note that when $k\rightarrow+\infty$, proper initialization in \eqref{l1tv-iteration} gives a global solution of the loss for all $\boldsymbol{y}_i$, and the network thus converges to the exact minimum of \eqref{l1tv-un}.
	However, when $k=L$ is fixed, or the initialization parameters are not chosen properly, the output of the network is generally not a solution of \eqref{l1tv-un}.
	Therefore, minimizing the empirical risk is actually to identify a parameter configuration that effectively alleviates the sub-optimality of the network with respect to \eqref{l1tv-un} across the input distribution utilized for training.
	In this way, the network is trained to acquire an algorithm that can effectively approximate the solution of \eqref{l1tv-un} for a specific category of signal distributions.
	Note that although the procedure can effectively expedite the resolution of the problem, the learned algorithm will only be effective for inputs $\boldsymbol{y}_i$ within the same input distribution as the training samples (i.e., in-distribution test).
	It may fail to deliver accurate approximations for samples that deviate significantly from the training set (i.e., out-of-distribution), in contrast to the iterative algorithm itself, which is guaranteed to converge regardless of the sample.
	
	In addition, the selection of penalty parameters holds great importance.
	Indeed, if $\lambda_1$ and $\lambda_2$ are big enough such that $\lambda_1+2\lambda_2\geq\lambda_{\rm max}$ (see Lemma \ref{l1tv-lem4.3} below), then the model \eqref{l1tv-un} has only  $\boldsymbol{0}$ as its solution.
	Therefore, it is crucial to select suitable penalty parameters for the model \eqref{l1tv-un}.
	Fortunately, unrolled networks naturally inherit prior structures and domain knowledge rather than learn this information from intensive training data \cite{monga2021algorithm}, and we  can initialize the penalty parameter pair $(\lambda_1,\lambda_2)$ based on the algorithm itself. 	
	\begin{lemma}\label{l1tv-lem4.3}
		Denote $\lambda_{\rm max}$ as the smallest value of the penalty parameter such that $\boldsymbol{0}$ is a solution. The value $\lambda_{\rm max}$ of the model \eqref{l1tv-un} is
		$\lambda_{\rm max}=\|\boldsymbol{A}^{\top}\boldsymbol{y}\|_\infty$.
	\end{lemma}
	\begin{proof}
		The first order conditions for problem \eqref{l1tv-un} reads: for $\boldsymbol{0}\in\mathbb{R}^n$
		\begin{equation*}
			\boldsymbol{0}\in\partial F(\boldsymbol{0}) = \boldsymbol{A}^{\top}\boldsymbol{y}+\lambda_1\partial\|\cdot\|_1(\boldsymbol{0})+\lambda_2\boldsymbol{D}^{\top}\|\cdot\|_1(\boldsymbol{D}\boldsymbol{0}).
		\end{equation*}
		Upon expanding the equality, we have
		\begin{equation*}
			\boldsymbol{a}_j^{\top}\boldsymbol{y}\in\begin{cases}
				[-(\lambda_1+\lambda_2),(\lambda_1+\lambda_2)],& j=1,\\
				[-(\lambda_1+2\lambda_2),(\lambda_1+2\lambda_2)],&\ j\in\{2,3,\cdots,n\}.
			\end{cases}
		\end{equation*}
		Thus, $\boldsymbol{0}$ is a solution of the model \eqref{l1tv-un} for all $\lambda_1+2\lambda_2\geq\lambda_{\rm max}$, with $\lambda_{\rm max}=\|\boldsymbol{A}^{\top}\boldsymbol{y}\|_\infty$.
	\end{proof}
	
	Note that LPGM-ISTA is designed to efficiently solve the regularized model \eqref{l1tv-un} with fixed $\lambda_1$ and $\lambda_2$, by unrolling PGM-ISTA for a fixed number of iterations. It is interesting to learn the parameters $\lambda_1$ and $\lambda_2$ simultaneously. For example, Scetbon et al \cite{Scetbon:2021} proposed a multi-layer perceptron (MLP) network for learning regularization parameters in the K-SVD algorithm. This idea may be useful for learning the parameters $\lambda_1$ and $\lambda_2$ in PGM-ISTA. We plan to explore the issue in future works.
	
	\section{Numerical experiments and discussions}\label{l1tv-simulation}
	
	In this section, we carry out some numerical experiments on the sampling number obtained in Section \ref{l1tv-section3} and the influence of parameters $(u,t)$ in PGM-ISTA, and the application of LPGM-ISTA to recover ECG signals.
	All the experiments are conducted on a desktop with Intel Core i9-13900K processor (3.00GHz, 64GB RAM) and NVIDIA GeForce RTX 4080 (16GB RAM).
	
	\subsection{Synthetic experiments for sampling number}
	In this part, we conduct synthetic experiments for the sampling number in Section \ref{l1tv-section3}, including a comparison with previous works. We shall consider $\epsilon=0$, and study the bound \eqref{l1tv-85} for the sampling number in several situations. Then we can fix $\lambda_2=1$ such that there is only $\lambda_1$ to be tuned. Note that PGM-ISTA is designed for the regularized $\ell^1$-TV method \eqref{l1tv-un} with $\epsilon\geq0$, which involve two penalty parameters $\lambda_1$ and $\lambda_2$. Moreover, since the constrained model \eqref{l1tv-model} with $\epsilon =0$ involves only one parameter and can be solved exactly by alternating direction method of multipliers (ADMM), we implement the model \eqref{l1tv-model} by ADMM instead of PGM-ISTA in this part.
	
	To develop the ADMM for the model \eqref{l1tv-model}, we rewrite as
	\begin{align}
		\min_{\boldsymbol{x}\in\mathbb{R}^n} g(\boldsymbol{x}), \quad \textup{s.t.} \quad \boldsymbol{y}=\boldsymbol{A}\boldsymbol{x}.
	\end{align}
	Let $\boldsymbol{x}=\boldsymbol{z}$, where $\boldsymbol{z}$ is an auxiliary variable. Then we obtain
	\begin{align}\label{l1tv-80}
		\min_{\boldsymbol{x},\boldsymbol{z}\in\mathbb{R}^n} g(\boldsymbol{z}), \quad \textup{s.t.} \quad \left[\begin{array}{c}
			\boldsymbol{A}\\
			\boldsymbol{I}_n
		\end{array}
		\right]\boldsymbol{x}=\left[ \begin{array}{c}
			\boldsymbol{y}\\
			\boldsymbol{z}
		\end{array}
		\right].
	\end{align}
	The augmented Lagrangian function is given by
	\begin{align*}
		g(\boldsymbol{z}) + \boldsymbol{u}^\top(\boldsymbol{Ax}-\boldsymbol{y})+\frac{\mu}{2}\|\boldsymbol{Ax}-\boldsymbol{y}\|_2^2 + \boldsymbol{v}^\top(\boldsymbol{x}-\boldsymbol{z})+\frac{\mu}{2}\|\boldsymbol{x}-\boldsymbol{z}\|_2^2,
	\end{align*}
	where $\boldsymbol{u}$ and $\boldsymbol{v}$ are dual variables.
	Following the procedure of ADMM, we update the variables alternatingly.
	\begin{itemize}
		\item For the term involving $\boldsymbol{x}$
		\begin{align}\label{l1tv-81}
			\boldsymbol{x}^{k+1}=&\mathop{\arg\min}_{\boldsymbol{x}}\boldsymbol{u}^{k,\top}(\boldsymbol{Ax}-\boldsymbol{y}) +\tfrac{\mu^k}{2}\|\boldsymbol{Ax}-\boldsymbol{y}\|_2^2 + \boldsymbol{v}^{k,\top}(\boldsymbol{x}-\boldsymbol{z}^k)+\tfrac{\mu^k}{2}\|\boldsymbol{x}-\boldsymbol{z}^k\|_2^2\nonumber\\ =&(\boldsymbol{A}^\top\boldsymbol{A}+\boldsymbol{I})^{-1}(\boldsymbol{A}^\top\boldsymbol{y}+\boldsymbol{z}^{k}- \tfrac{1}{\mu^{k}}\boldsymbol{A}^\top\boldsymbol{u}^{k}-\tfrac{1}{\mu^{k}}\boldsymbol{v}^{k}).
		\end{align}
		\item For the term involving $\boldsymbol{z}$, by Lemma \ref{l1tv-lemma4}, we have
		\begin{align}\label{l1tv-82}
			\boldsymbol{z}=&\mathop{\arg\min}_{\boldsymbol{z}}g(\boldsymbol{z})+\tfrac{\mu^k}{2}\|\boldsymbol{x}^{k+1}- \boldsymbol{z}+\tfrac{\boldsymbol{v}^{k}}{\mu^{k}}\|_2^2
			=\mathcal{P}_{\lambda_1/\mu^{k}}^{\lambda_2/\mu^{k}}\Big(\boldsymbol{x}^{k+1}+\tfrac{\boldsymbol{v}^{k}}{\mu^{k}}\Big).
		\end{align}
	\end{itemize}
	We summarize the whole update procedure in Algorithm \ref{l1tv-admm}.
	Since this is a two block ADMM, the global convergence is guaranteed \cite{boyd2011distributed}.

	\begin{algorithm}[ht]
		\caption{ADMM for problem \eqref{l1tv-80}}
		\label{l1tv-admm}
		\begin{algorithmic}
			\STATE {\textbf{Input}: Sensing matrix $\boldsymbol{A}$, measurement $\boldsymbol{y}$.}\\
			\ 1: Initialization: $\boldsymbol{x}^0,\boldsymbol{z}^0$ being zero vectors, regularization parameters $\lambda_1$, $\lambda_2$, tolerate error $tol=10^{-8}$,  $\mu^0=10^{-3}$, $\rho=1.1$, $\mu_{max}=10^8$, and $k=0$.\\
			\textbf{while} not convergent \textbf{do} \\
			\ 2: Update $\boldsymbol{x}^{k+1}$ by \eqref{l1tv-81};\\
			\ 3: Update $\boldsymbol{z}^{k+1}$ by \eqref{l1tv-82};\\
			\ 4: $\boldsymbol{u}^{k+1} = \boldsymbol{u}^{k} + \mu^{k}(\boldsymbol{A}\boldsymbol{x}^{k+1} - \boldsymbol{y})$;\\
			\ 5: $\boldsymbol{v}^{k+1} = \boldsymbol{v}^{k} + \mu(\boldsymbol{x}^{k+1}-\boldsymbol{z}^{k+1})$;\\
			\ 6: $\mu^{k+1} = \min\{\rho\mu^k,\mu_{max}\}$; \\
			\ 7: Check the convergence condition\\
			\quad~$\|\boldsymbol{x}^{k+1}-\boldsymbol{x}^k\|_2/\max\{\|\boldsymbol{x}^k\|_2,1\} < tol$;\\
			\ 8: Update $k\leftarrow k+1$;\\
			\ 9: \textbf{end while}\\
			\textbf{Output:} $\hat{\boldsymbol{x}}=\boldsymbol{x}^{k+1}$.
		\end{algorithmic}
	\end{algorithm}

\begin{figure*}[hbt!]
		\renewcommand{\arraystretch}{0.5}
		\setlength\tabcolsep{0.5pt}
		\centering
		\begin{tabular}{ccccccc}
			\centering
			\includegraphics[width=66mm, height = 50mm]{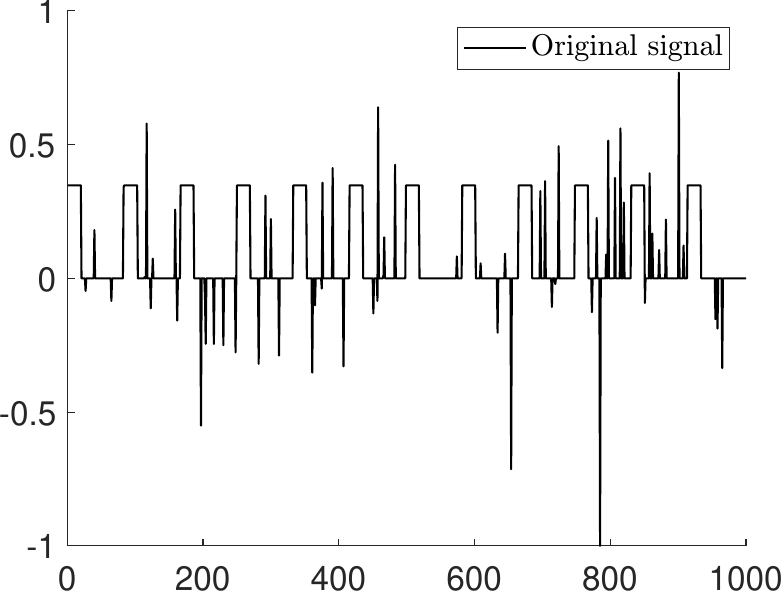}&\qquad
			\includegraphics[width=66mm, height = 50mm]{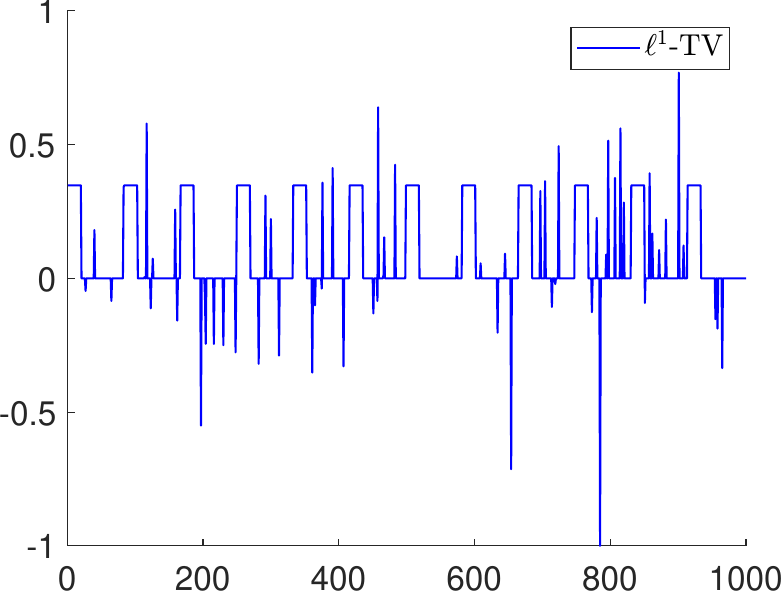}
		\end{tabular}
		\caption{A synthetic signal and its recovered one by the $\ell^1$-TV method.}
		\label{syn_syn_0}
	\end{figure*}

	Next, we carry out simulations to complement the theory in Section \ref{l1tv-section3}.
	Since the $\ell^1$-TV model enforces sparsity both in signal itself and the differences between
	consecutive components of the signal, we conduct experiments on signals that are locally constant and change in jumps, following the setting in \cite{ye2011split}. Specifically, we generate the synthetic signal ${\boldsymbol{x}}^\ast\in\mathbb{R}^n$ in four steps: (i) Set $\|{\boldsymbol{x}}^\ast\|_0=s_r<n$; (ii) To ensure the local-smooth property of $\boldsymbol{x}^\ast$, we generate $b$ ($b\in\mathbb{N}$) blocks randomly with block size 10, and elements in each block have the same amplitude, sampled from the standard normal distribution $\mathcal{N}(0,1)$;
	(iii) we generate $s_r-10b$ locations with their amplitudes following $\mathcal{N}(0,1)$;
	(iv) we normalize the signal $\boldsymbol{x}^\ast$ to the interval $[-1,1]$ by $\boldsymbol{x}^\ast/\max|\boldsymbol{x}^\ast|$.
	The Gaussian sampling strategy on $\boldsymbol{A}$ stated in Theorem \ref{l1tv-theorem3.2} is adopted.
	To illustrate the recovery accuracy, we show in Fig. \ref{syn_syn_0} one synthetic signal with $n=1000$, $s_r=300$ and $b=12$, and its recovered signal by $\ell^1$-TV method under $\epsilon=0$ with a sampling ratio of 0.5.
	
	\begin{figure*}[hbt!]
		\renewcommand{\arraystretch}{0.5}
		\setlength\tabcolsep{0.5pt}
		\centering
		\begin{tabular}{ccccccc}
			\centering
			\includegraphics[width=40mm, height = 40mm]{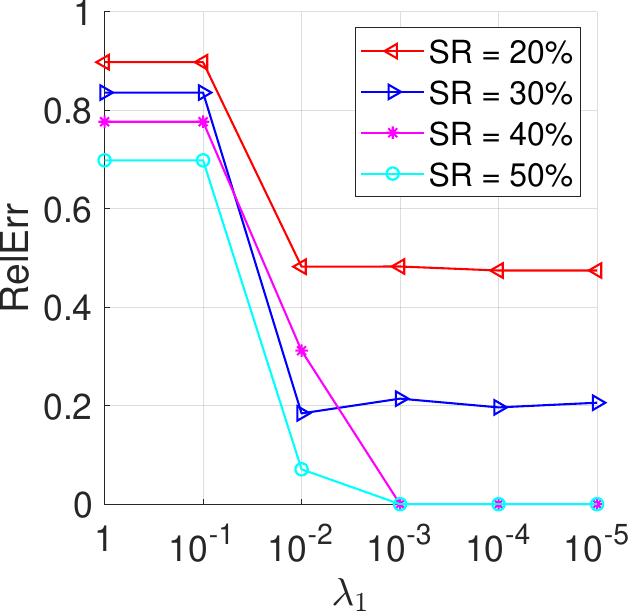}&
			\includegraphics[width=40mm, height = 40mm]{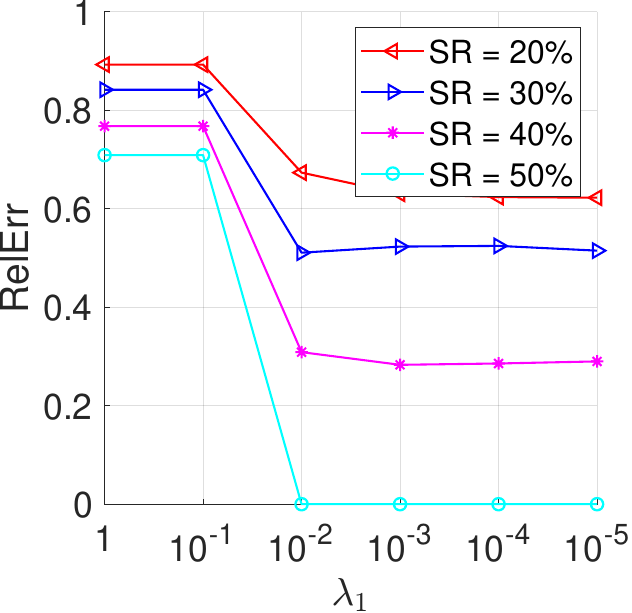}&
			\includegraphics[width=40mm, height = 40mm]{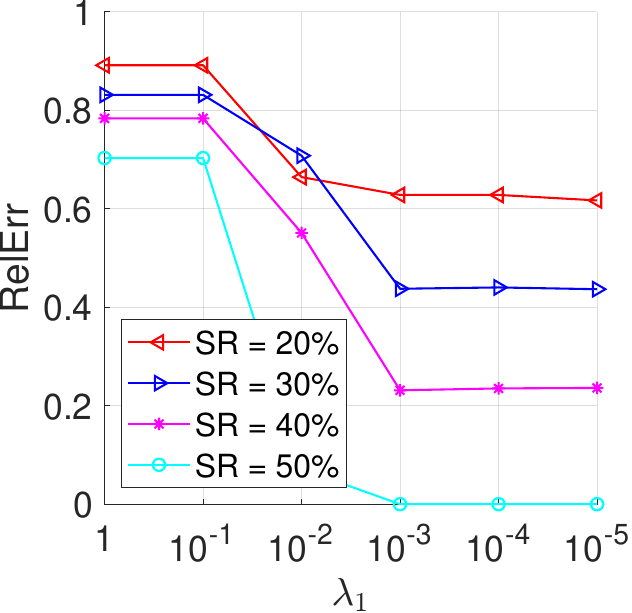}&
			\includegraphics[width=40mm, height = 40mm]{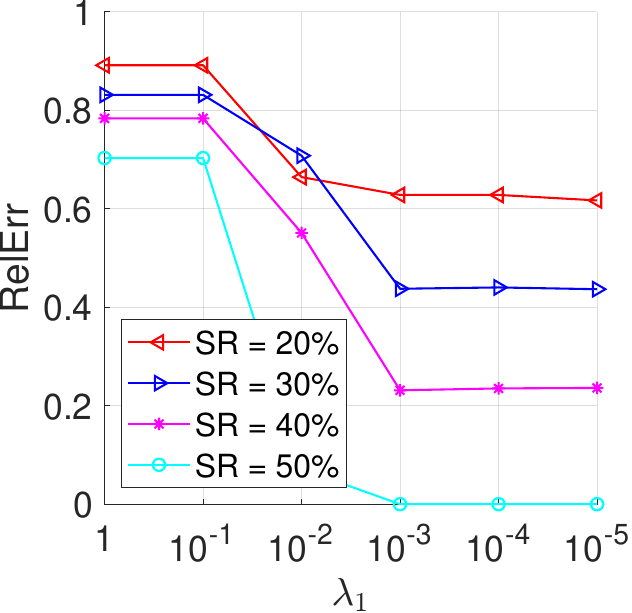}\\
			\scriptsize $n=500,s_r=100,s_g=44$&
			\scriptsize $n=500,s_r=150,s_g=66$&
			\scriptsize $n=1000,s_r=200,s_g=91$&
			\scriptsize $n=1000,s_r=300,s_g=137$
		\end{tabular}
		\caption{Selections of $\lambda_1$ for fixed $\lambda_2=1$ under different situations.}
		\label{syn_lambda1}
	\end{figure*}
	
	\begin{figure*}[hbt!]
		\renewcommand{\arraystretch}{0.5}
		\setlength\tabcolsep{0.5pt}
		\centering
		\begin{tabular}{ccccccc}
			\centering
			\includegraphics[width=40mm, height = 40mm]{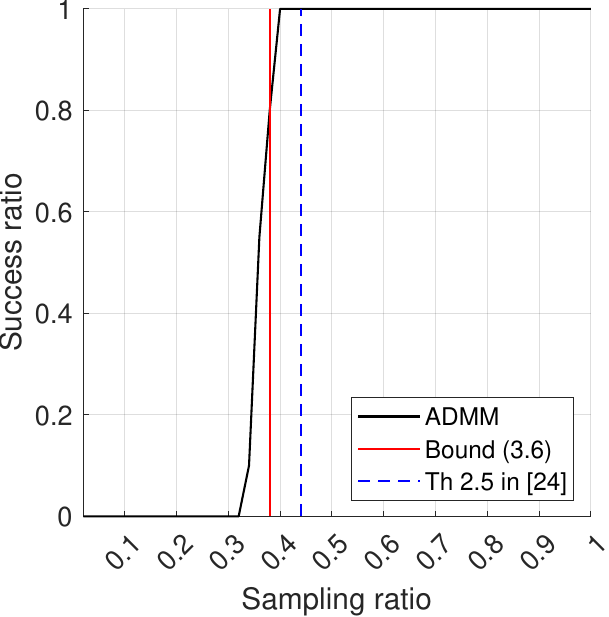}&
			\includegraphics[width=40mm, height = 40mm]{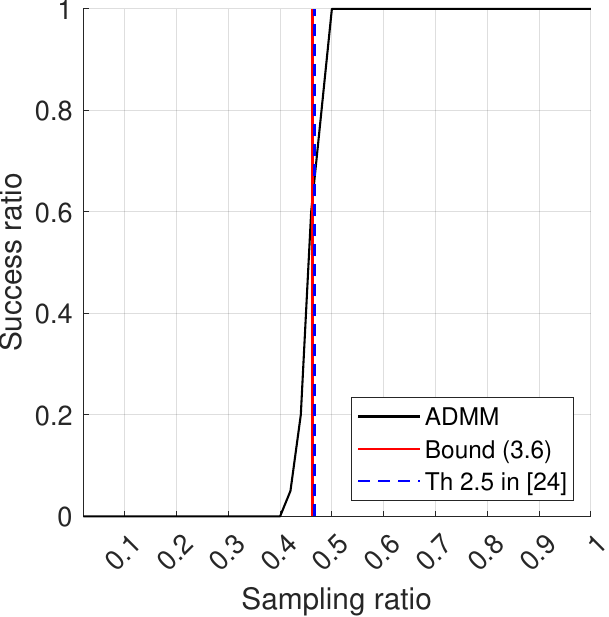}&
			\includegraphics[width=40mm, height = 40mm]{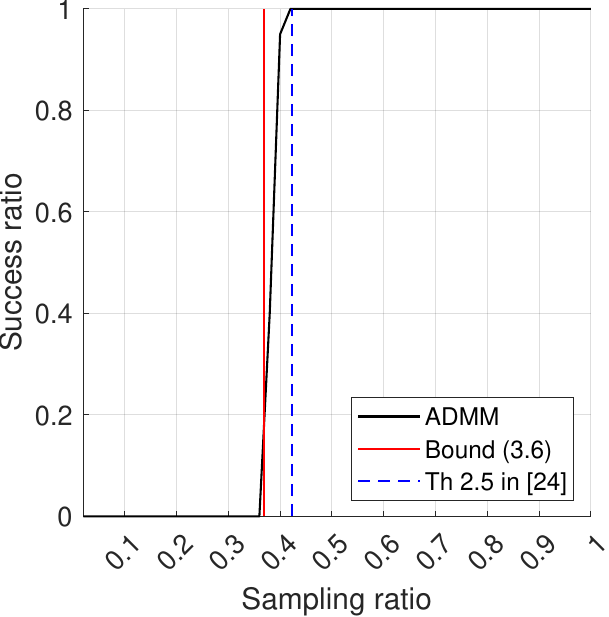}&
			\includegraphics[width=40mm, height = 40mm]{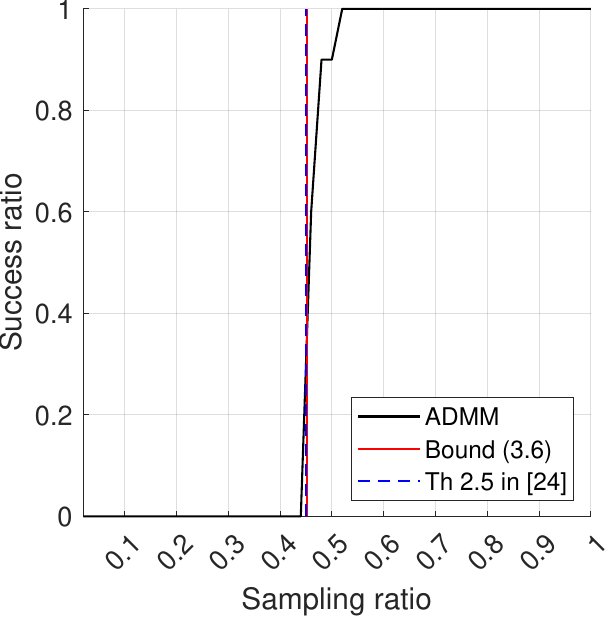}\\
			\includegraphics[width=40mm, height = 40mm]{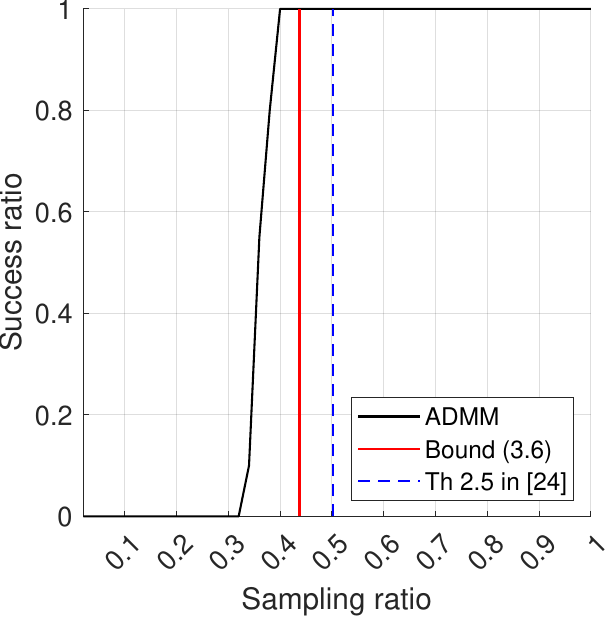}&
			\includegraphics[width=40mm, height = 40mm]{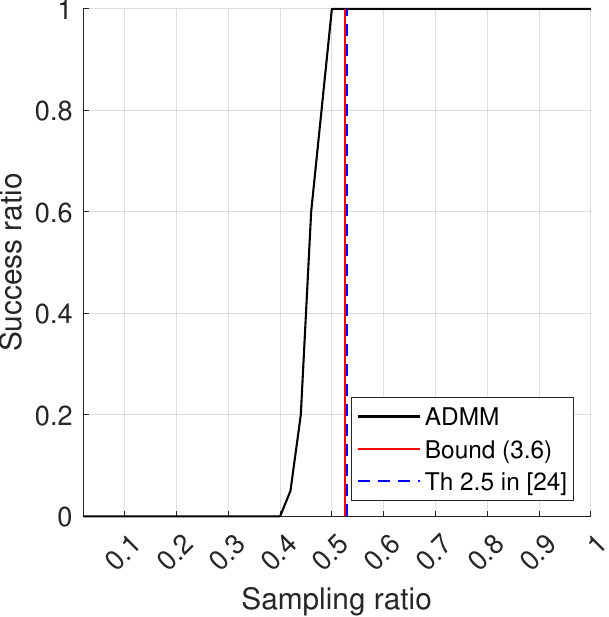}&
			\includegraphics[width=40mm, height = 40mm]{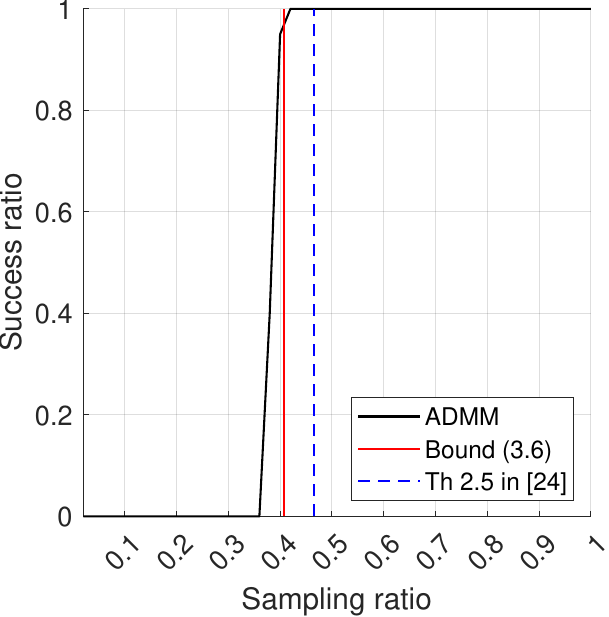}&
			\includegraphics[width=40mm, height = 40mm]{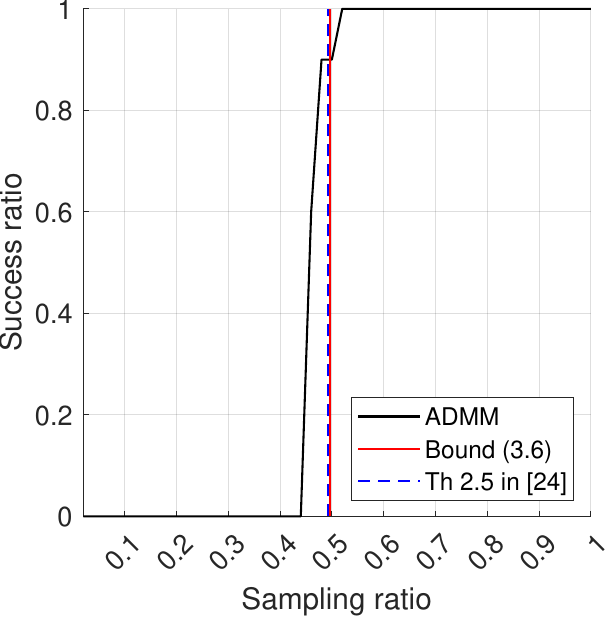}\\
			\scriptsize $n=500,s_r=100,s_g=44$&
			\scriptsize $n=500,s_r=150,s_g=66$&
			\scriptsize $n=1000,s_r=200,s_g=91$&
			\scriptsize $n=1000,s_r=300,s_g=137$
		\end{tabular}
		\caption{The theoretical bound \eqref{l1tv-85} and the bound from Theorem 2.5 of \cite{genzel2021l1} at different sparsity levels, with $t=1$ (top row) and  $t=2$ (bottom row).}
		\label{l1tv-phase}
	\end{figure*}
	
	We first examine the influences of the parameters $\lambda_1$ and $\lambda_2$ on the recovery performance.
	We simply set $\lambda_2=1$ in the model \eqref{l1tv-model}, and thus we only test the influence of $\lambda_1$.
	To this end, we choose $\lambda_1$ from the set $\{1,10^{-1},10^{-2},10^{-3},10^{-4},10^{-5}\}$ by cross-validation.
	Moreover, we consider different $n$, $s_r$ and $s_g$, and plot the relative errors versus $\lambda_1$ in Fig. \ref{syn_lambda1}.
	For a recovered signal $\hat{{\boldsymbol{x}}}$ with the clean signal ${\boldsymbol{x}}^\ast$, we define the relative error (RelErr) as $	 \textup{RelErr}(\hat{\boldsymbol{x}},\boldsymbol{x}^\ast) = \|\hat{{\boldsymbol{x}}} - {\boldsymbol{x}}^\ast\|_2/\|{\boldsymbol{x}}^\ast\|_2$.
	It is observed that $\lambda_1=10^{-3}$ is a good choice in all cases.
	Further, for different signal lengths and sparsity levels, we examine the influence of the sampling number $m$ on the recovery performance.
	To this end, for each sampling ratio $m/n$, we repeat the experiments 50 times and declare success for the recovered signal $\hat{{\boldsymbol{x}}}$ if $\|{\boldsymbol{x}^\ast}-\hat{{\boldsymbol{x}}}\|_2<10^{-3}$.
	The curves in Fig. \ref{l1tv-phase} show the success ratio with respect to the sampling ratio, along with the bound \eqref{l1tv-85} and the bound stated in Theorem 2.5 of \cite{genzel2021l1}, representing the required number of scaled form samples for successful recovery.
	Compared with the bound due to \cite{genzel2021l1} under the condition of $t=1$ (with a probability $1-e^{-t^2/2}\approx0.3935$) and $t=2$ (with a probability $1-e^{-t^2/2}\approx0.8647$), the bound \eqref{l1tv-85} can more faithfully estimate the sampling number $m$ needed for successful recovery, which appears tighter in the simulations, especially for relatively low levels of regular sparsity and gradient sparsity.
	
	\subsection{Parameters influence for PGM-ISTA}
	
	Now we study the influence of the parameter pair $(u,t)$ on the convergence of PGM-ISTA.
	We choose one ECG signal from the MIT-BIH Arrhythmia Database \cite{moody2001the} as
	the signal to be recovered. The ECG signal is windowed to 1024 samples, cf. Fig. \ref{l1tv-intro-show}.
	The sensing matrix $\boldsymbol{A}$ is Gaussian, and the measurements are contaminated by Gaussian noise with variance $0.01$ and a sampling ratio 0.5.
	
	\begin{figure*}[hbt!]
		\renewcommand{\arraystretch}{0.5}
		\setlength\tabcolsep{0.5pt}
		\centering
		\begin{tabular}{ccccccc}
			\centering
			\includegraphics[width=40mm, height = 40mm]{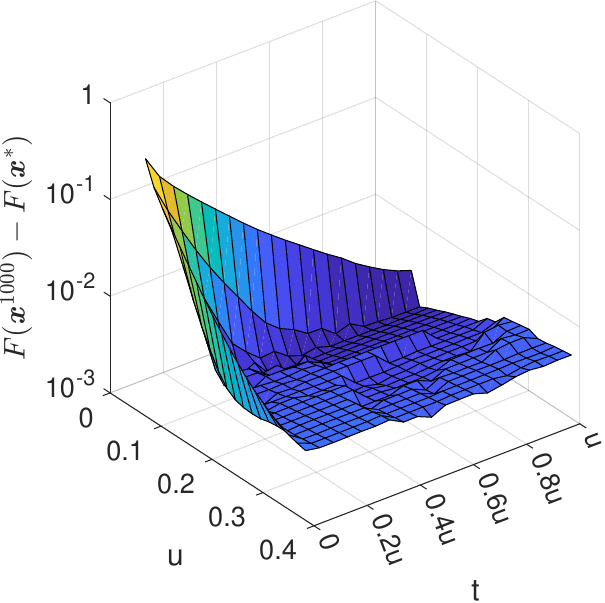}&
			\includegraphics[width=40mm, height = 40mm]{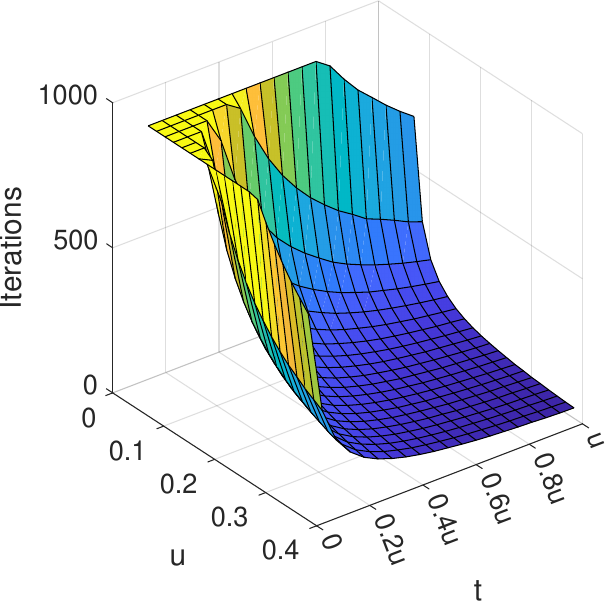}&
			\includegraphics[width=40mm, height = 40mm]{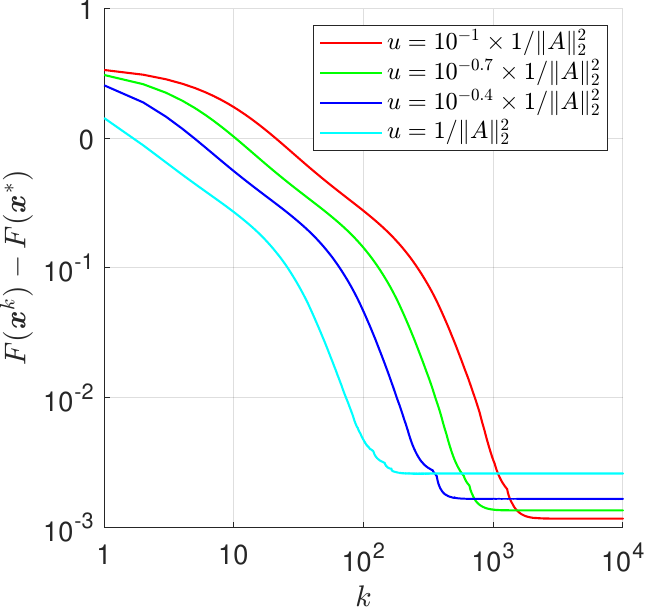}&
			\includegraphics[width=40mm, height = 40mm]{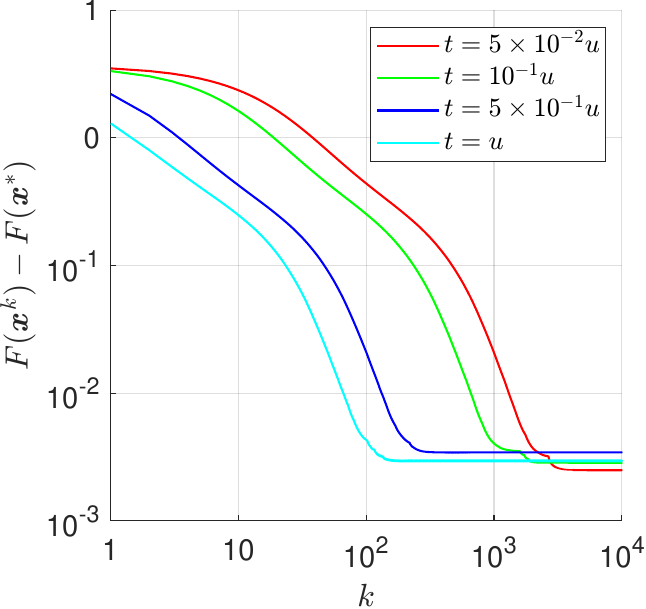}\\
			\scriptsize (a)&
			\scriptsize (b)&
			\scriptsize (c)&
			\scriptsize (d)
		\end{tabular}
		\caption{The convergence of PGM-ISTA for fixed $\lambda_1 = \lambda_2 = 0.01$, (a) shows values of $F(\boldsymbol{x}^{1000}(u,t)) - F(\boldsymbol{x}^\ast)$; (b) shows the minimum iteration $k$ such that $F(\boldsymbol{x}^{k}(u,t)) - F(\boldsymbol{x}^\ast)<0.01$; (c) shows values of $F(\boldsymbol{x}^{k}(u,t)) - F(\boldsymbol{x}^\ast)$ with multiple choices of $u$ and $t=0.9u$; (d) shows values of $F(\boldsymbol{x}^{k}(u,t)) - F(\boldsymbol{x}^\ast)$ for multiple choices of $t$ and $u=1/{\|\boldsymbol{A}\|_2^2}$.}
		\label{l1tv-alg-ut}
	\end{figure*}
	
	The PGM-ISTA algorithm and its learned counterpart include one parameter pair $(u,t)$ and the model \eqref{l1tv-un} contains two
	penalty parameters $\lambda_1$ and $\lambda_2$. The first experiment aims to verify the existence
	of a parameter pair $(u,t)$ by cross-validation for fixed penalty parameters.
	We take $\lambda_1 = \lambda_2 = 0.01$, and conduct an experiment with both $u$ and $t$ being
	in an increasing order and equi-spaced for $u\in(0,2/\|\boldsymbol{A}\|_2^2)$ and $t\in(0,u]$
	by Theorem \ref{l1tv-thm2}, respectively. For each pair $(u,t)$, we set the maximum number of
	iterations to 1000, and the results are shown in Fig. \ref{l1tv-alg-ut}.
	Fig. \ref{l1tv-alg-ut} (a) shows that there are many pairs $(u,t)$ for which PGM-ISTA achieves the goal
	$F({\boldsymbol{x}}^k) - F({\boldsymbol{x}}^\ast) < 0.01$ within 1000 iterations.
	Thus, there is a large region of $(u,t)$ providing convergent sequences for
	both PGM-ISTA if we declare convergence under the criterion $F({\boldsymbol{x}}^k) -
	F({\boldsymbol{x}}^\ast) < 0.01$ for some $k\in\mathbb{N}_+$. Fig. \ref{l1tv-alg-ut} (b)
	shows the number of PGM-ISTA iterations needed to achieve $F({\boldsymbol{x}}^k) -
	F({\boldsymbol{x}}^\ast) < 0.01$. Clearly, bigger values of pair $(u,t)$ may result in faster convergence.
	The numerical results contrast slightly the theoretical prediction of Theorem \ref{l1tv-thm4} which suggests that a smaller $u$ may give a smaller error in the function value for PGM-ISTA.
	This is attributed to the fact that we only run PGM-ISTA for a maximum of 1000 iterations, which seems not enough occasionally.
	Further experiments indicate  that a small $u$ may lead to slower convergence. Hence, we have increased the number of iterations for several choices of $u$, and the results are shown in Fig. \ref{l1tv-alg-ut} (c).
	Then indeed a smaller $u$ leads to a smaller function error, which agrees well with the
	theoretical results from Theorem \ref{l1tv-thm4}, but a smaller $u$ does require more iterations to reach the desired convergence.
	
	Next, we explore how the parameter $t$ affects the convergence rate in terms of the function value.
	We conduct an experiment for fixed $u=1/\|\boldsymbol{A}\|_2^2$ and different $t$, and Fig. \ref{l1tv-alg-ut} (d) shows the results.
	Obviously, the influence of $t$ about the function value is slight, but indeed a bigger $t$ may lead to faster convergence.
	
	Finally, we compare the performance of the $\ell^1$, TV and $\ell^1$-TV methods when the observation data are contaminated by additive i.i.d. Gaussian noise (with a standard deviation $\sigma$).
	We choose $\sigma$ from the set $\{0,0.01,0.1\}$ and the signal in Fig. \ref{syn_syn_0} as the ground truth, and summarize the relative errors of the recovered signals by the three methods in Table \ref{l1tv-tab2}. It is observed that $\ell^1$-TV can obtain much smaller relative errors in all cases, showing its superiority when the exact signal exhibits both signal sparsity and gradient sparsity. Also we plot the exact signal and the recovered ones with $\sigma=0$ in Fig. \ref{syn_syn_1}, which again confirms that $\ell^1$-TV performs better than the $\ell^1$ or TV methods. Note that the algorithm ADMM cannot be used when $\sigma\neq0$.
	
	\begin{table}[hbt!]
		\setlength{\tabcolsep}{1.8mm}
		\centering
		\fontsize{7}{8}\selectfont
		\begin{threeparttable}
			
			\caption{Relative errors of three methods under $\sigma\in\{0,0.01,0.1\}$.\label{l1tv-tab2}}
			\begin{tabular}{ccccccccccccccc}
				\toprule[1pt]
				Methods&$\sigma=0$&$\sigma=0.01$&$\sigma=0.1$\cr
				\hline
				$\ell^1$&0.7000&0.7488&0.7516&\cr
				TV&0.1025&0.1224&0.1296&\cr
				$\ell^1$-TV&0.0027&0.0088&0.0793\cr
				\bottomrule[1pt]
			\end{tabular}
		\end{threeparttable}
		
	\end{table}
	
	\begin{figure*}[hbt!]
		\renewcommand{\arraystretch}{0.5}
		\setlength\tabcolsep{0.5pt}
		\centering
		\begin{tabular}{ccccccc}
			\centering
			\includegraphics[width=40mm, height = 31mm]{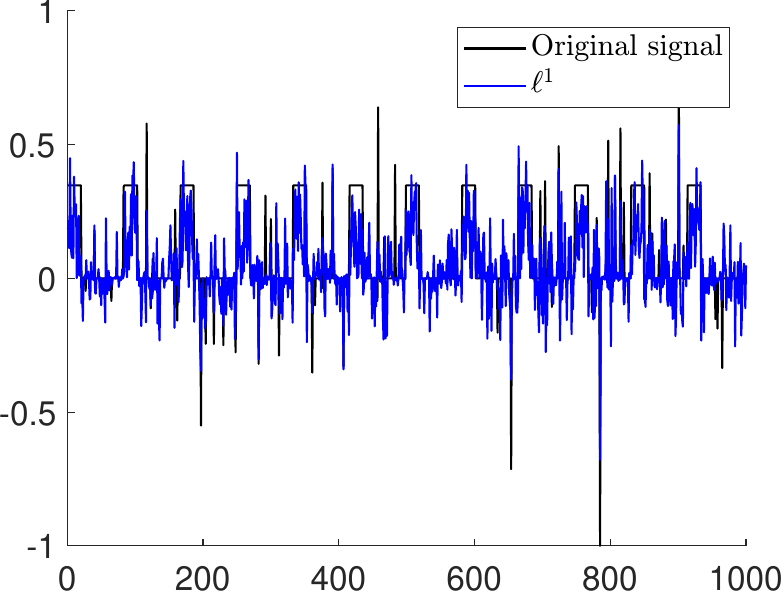}&
			\includegraphics[width=40mm, height = 31mm]{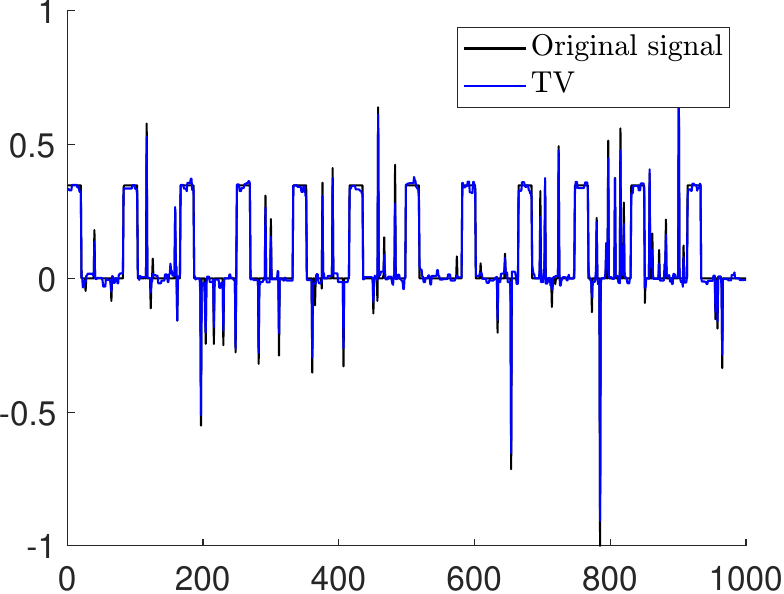}&
			\includegraphics[width=40mm, height = 31mm]{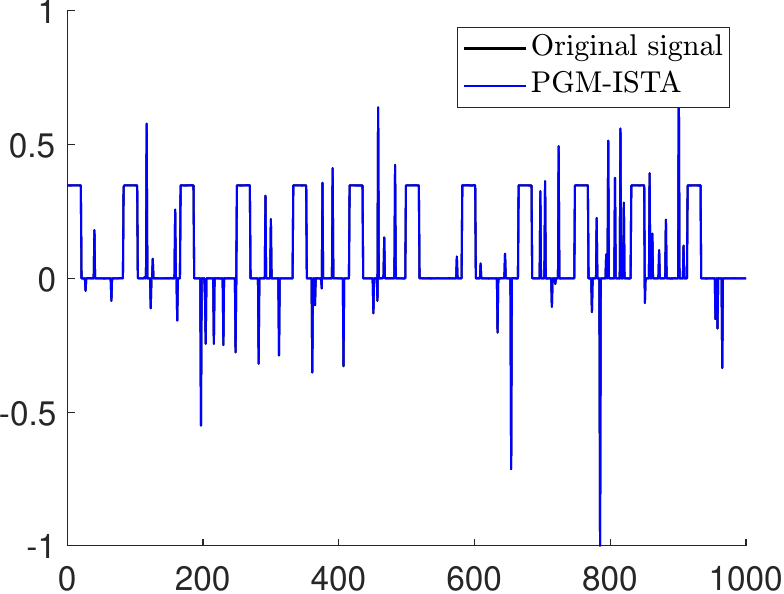}&
			\includegraphics[width=40mm, height = 31mm]{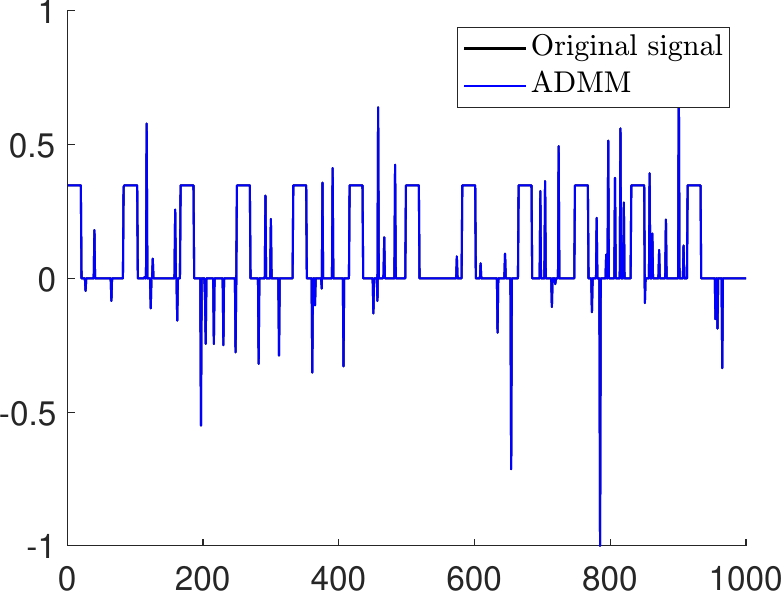}\\
			\scriptsize (a)&
			\scriptsize (b)&
			\scriptsize (c)&
			\scriptsize (d)
		\end{tabular}
		\caption{The exact and recovered signals by three methods with a sampling ratio 0.5 and $\sigma=0$. 	\label{syn_syn_1}
		}
	\end{figure*}
	
	\subsection{LPGM-ISTA for recovering ECG signals}
	
	Last we conduct experiments to recover ECG signals, to illustrate the effectiveness of the $\ell^1$-TV model \eqref{l1tv-un} on real-world signals, and compare LPGM-ISTA with traditional iterative algorithms.
	All experiments of neural networks are performed in Python using PyTorch.
	
	We choose 2372 ECG signals $\{\boldsymbol{x}_i\}_{i=1}^{2372}$ with length windowed to 256 from the MIT-BIH Arrhythmia Database \cite{moody2001the}.
	The observations are generated by a linear transform as $\boldsymbol{y}_i=\boldsymbol{A}{\boldsymbol{x}}_i$, where $\boldsymbol{A}\in\mathbb{R}^{128\times256}$ with entries following i.i.d. $\mathcal{N}(0,1)$.
	In this way, we obtain the pair of observed signals and its labels $\{(\boldsymbol{y}_i,{\boldsymbol{x}}_i)\}_{i=1}^{2366}$.
	Next, we choose the first 1900 samples for training, and the remaining 472 samples for testing.
	The regularization parameters $\lambda_1$ and $\lambda_2$ are set according to cross-validation, which are verified to satisfy the condition $\lambda_1+2\lambda_2<\lambda_{\rm max}$.
	For a typical signal from the test set, we compare the performance of LPGM-ISTA and several existing algorithms for solving problem \eqref{l1tv-un} including the proposed PGM-ISTA, LADMM \cite{li2014linearized}, and S-FISTA \cite{beck2012smoothing}.
	
	We report the relative errors and CPU time in Table \ref{l1tv-alg-len} with iteration / layer number ranging from 2 to 10 with a step size of 2, and also the results for LADMM, S-FISTA and PGM-ISTA with the iteration number being 500 and 1000. Clearly, LPGM-ISTA far outperforms the three iterative algorithms including LADMM, S-FISTA and the proposed PGM-ISTA in terms of both relative errors and CPU time. The relative errors of LPGM-ISTA show a decreasing trend with the increase of layer / iteration number, and moreover very few layers of the learned algorithm (LPGM-ISTA) are sufficient to achieve satisfactory relative errors compared with the original algorithm, whereas the iterative algorithms lose their effectiveness in terms of the required number of iterations. The exact signal and the recovered ones by the four algorithms with iteration / layer number being 2 are shown in Fig. \ref{l1tv-ecg-rev}.
	Most remarkably, LPGM-ISTA can already provide a satisfactory recovery, but LADMM, S-FISTA and PGM-ISTA barely recover any useful information.
	Thus, LPGM-ISTA does provide a fast solver for problem \eqref{l1tv-un}.
	Here the signal length is only 256, which is relatively short. It is expected that the proposed fast learned solver may exhibit greater superiority when dealing with even longer signals.
	
	\begin{table}[hbt!]
		\setlength{\tabcolsep}{1.8mm}
		\centering
		\fontsize{7}{8}\selectfont
		\begin{threeparttable}
			\caption{Comparisons of the RelErrs/time (in seconds) of LADMM, S-FISTA, PGM-ISTA and LPGM-ISTA.}
			\label{l1tv-alg-len}
			\begin{tabular}{cccccccc}
				\toprule[1pt]
				Iterations/layers&2&4&6&8&10&500&1000\cr
				\midrule[0.8pt]
				LADMM \cite{li2014linearized}&0.8277/-&0.7842/-&0.7615/-&0.7464/-&0.7345/-&0.1386/0.152&0.0391/0.295\cr
				S-FISTA \cite{beck2012smoothing}&0.8282/-&0.7720/-&0.7427/-&0.7266/-&0.7155/-&0.0712/0.184&0.0502/0.365\cr
				PGM-ISTA&0.8255/-&0.7801/-&0.7563/-&0.7413/-&0.7306/-&0.1435/0.204&0.0446/0.364\cr
				LPGM-ISTA&\bf0.0560/0.063&\bf0.0503/0.073&\bf0.0435/0.081&\bf0.0413/0.085&\bf0.0391/0.095&-/-&-/-\cr
				\bottomrule[1pt]
				\vspace{-0.8cm}
			\end{tabular}
		\end{threeparttable} 
	\end{table}	
	
	\begin{figure}[hbt!]
		\renewcommand{\arraystretch}{0.5}
		\setlength\tabcolsep{0.5pt}
		\centering
		\begin{tabular}{ccccccc}
			\centering
			\includegraphics[width=40mm, height = 31mm]{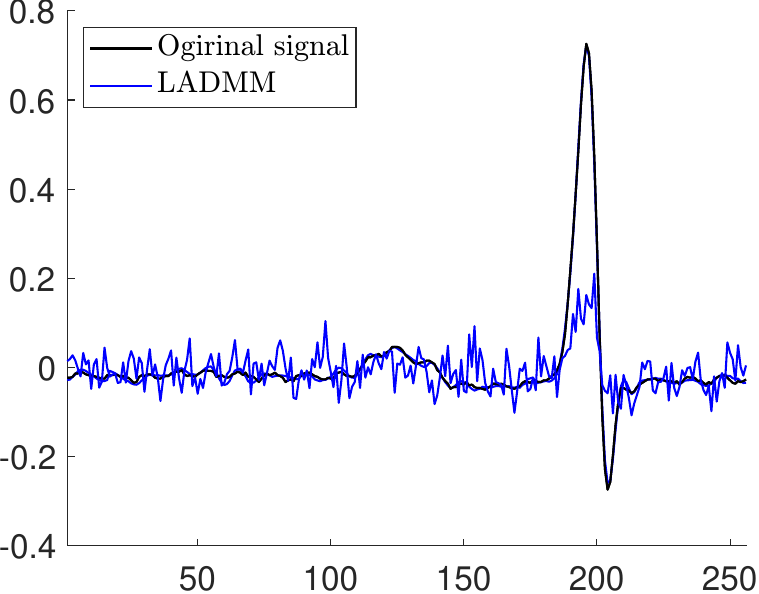}&
			\includegraphics[width=40mm, height = 31mm]{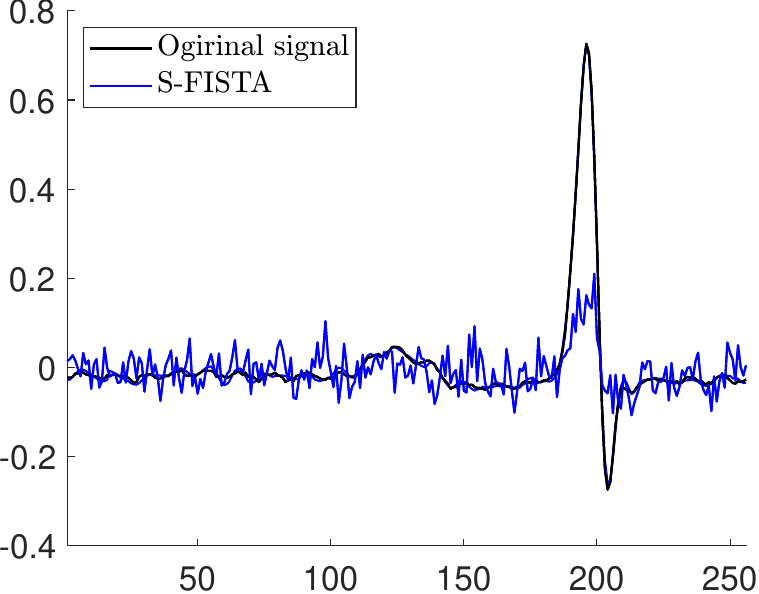}&
			\includegraphics[width=40mm, height = 31mm]{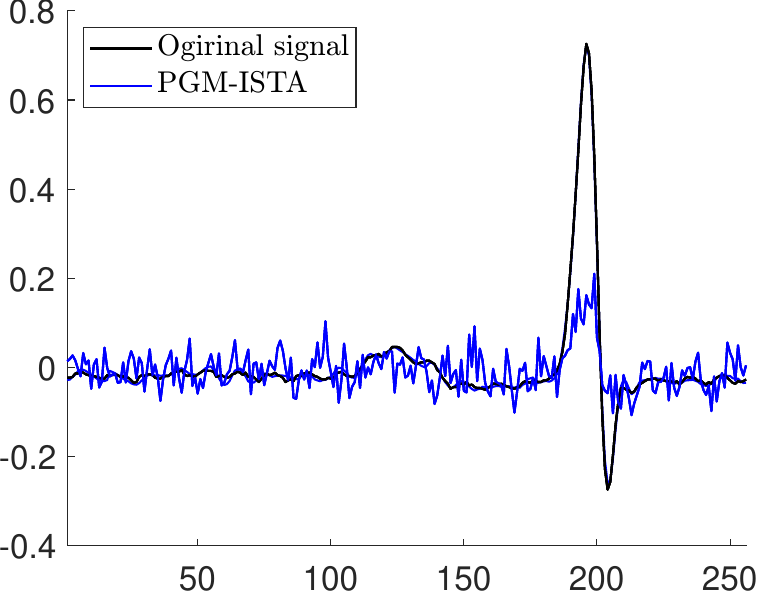}&
			\includegraphics[width=40mm, height = 31mm]{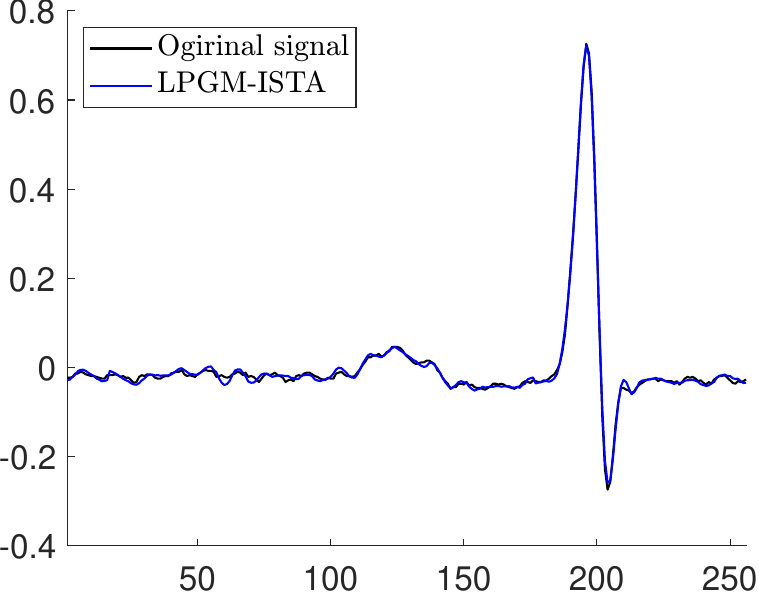}\\
			\scriptsize (a)&
			\scriptsize (b)&
			\scriptsize (c)&
			\scriptsize (d)
		\end{tabular}
		\caption{The exact signal and the recovered ones by LADMM, S-FISTA, PGM-ISTA and LPGM-ISTA with iterations/layers being 2. \label{l1tv-ecg-rev}}
	\end{figure}

	\section{Technical proofs}\label{l1tv-proof}
	
	In this section we give all the technical proofs.
	
	\subsection{Proof of Theorem \ref{l1tv-theorem3.1}}\label{l1tv-proof1}
	
	The proof is partly inspired from \cite{daei2018sample}. However, we need to make several essential changes in order to adapt the results to the $\ell^1$-TV case. Technically, we select a vector from the subdifferential of the objective function in \eqref{l1tv-model} that yields an upper bound on its statistical dimension. First we give two useful lemmas.
	\begin{lemma}[{\cite[Grothendieck's identity]{vershynin2019high}}]\label{l1tv-lemma2}
		For any fixed vectors $\boldsymbol{u}$, $\boldsymbol{v}\in \mathbb{S}^{n-1}$, we have
		\begin{equation} \operatorname{E}_{\boldsymbol{g}\sim \mathcal{N}(\boldsymbol{0},\boldsymbol{I}_n)}[\operatorname{sign}\langle\boldsymbol{g},\boldsymbol{u}\rangle\operatorname{sign} \langle\boldsymbol{g},\boldsymbol{v}\rangle] =\tfrac{2}{\pi}\arcsin\langle\boldsymbol{u},\boldsymbol{v}\rangle.
		\end{equation}
	\end{lemma}
	
	\begin{lemma}\label{l1tv-lemma9}
		Let $\boldsymbol{D}=[\boldsymbol{d}_1,\boldsymbol{d}_2,\cdots,\boldsymbol{d}_{n-1}]^\top$. Then for any $i\in[n-1]$, we have
		\begin{align*} \operatorname{E}_{\boldsymbol{g}\sim \mathcal{N}(\boldsymbol{0},\boldsymbol{I}_n)}[(\boldsymbol{d}_i^{\top}\operatorname{sign}(\boldsymbol{g}) \times\operatorname{sign}(\boldsymbol{d}_i^{\top}\boldsymbol{g}))]=0.
		\end{align*}
	\end{lemma}
	\begin{proof}
		The result follows directly by computing the expectation.
	\end{proof}
	
	Now we can state the proof of Theorem \ref{l1tv-theorem3.1}.
	Let $s_1$ and $s_2$ be the number of adjacent pairs in $\mathbb{S}_G$ and $\mathbb{S}_G^c$, respectively, with $\mathbb{S}_G$ denoting $\mathbb{S}_G({\boldsymbol{x}})$. That is,
	\begin{align*}
		s_1=\left|\{i\in\{2,3,\cdots,n-1\}:i\in \mathbb{S}_G,i-1\in \mathbb{S}_G\}\right|,\\
		s_2=\left|\{i\in\{2,3,\cdots,n-1\}:i\in \mathbb{S}_G^c,i-1\in \mathbb{S}_G^c\}\right|.
	\end{align*}
	Then $s_1\leq|\mathbb{S}_G|-1$ and $s_2\leq|\mathbb{S}_G^c|-1$.
	Since $\partial\|\cdot\|_1(\cdot)$ is a compact set, the following two optimization problems
	\begin{equation*}
		\boldsymbol{z}_1=\mathop{\arg\max}\limits_{\boldsymbol{z}\in\partial\|\cdot\|_1({\boldsymbol{x}})}\langle\boldsymbol{g},\boldsymbol{z}\rangle \quad \text{and} \quad \boldsymbol{z}_2=\mathop{\arg\max}\limits_{\boldsymbol{z}\in\partial\|\cdot\|_1(\boldsymbol{D}{\boldsymbol{x}})} \langle\boldsymbol{g},\boldsymbol{D}^{\top}\boldsymbol{z}\rangle
	\end{equation*}
	are well defined. In addition, we have
	\begin{align*}		 &\boldsymbol{z}_1=\mathop{\arg\max}\limits_{\boldsymbol{z}\in\partial\|\cdot\|_1({\boldsymbol{x}})}\langle\boldsymbol{g},\boldsymbol{z}\rangle=\mathop{\arg\max}\limits_{\|\boldsymbol{z}\|_\infty\leq 1}\langle\boldsymbol{g},\operatorname{sign}({\boldsymbol{x}})_{\mathbb{S}_R}+\boldsymbol{z}_{\mathbb{S}_R^c}\rangle=\operatorname{sign}({\boldsymbol{x}})_{\mathbb{S}_R}+\operatorname{sign}(\boldsymbol{g})_{\mathbb{S}_R^c}\\
		&\boldsymbol{z}_2=\mathop{\arg\max}\limits_{\boldsymbol{z}\in\partial\|\cdot\|_1({\boldsymbol{Dx}})}\langle\boldsymbol{g},\boldsymbol{D}^{\top}\boldsymbol{z}\rangle=\mathop{\arg\max}\limits_{\|\boldsymbol{z}\|_\infty\leq 1}\langle \boldsymbol{Dg},\operatorname{sign}({\boldsymbol{Dx}})_{\mathbb{S}_G}+\boldsymbol{z}_{\mathbb{S}_G^c}\rangle=\operatorname{sign}({\boldsymbol{Dx}})_{\mathbb{S}_G}+\operatorname{sign}(\boldsymbol{Dg})_{\mathbb{S}_G^c}.
	\end{align*}
	Consequently,
	\begin{align}
		\operatorname{dist}^2(\boldsymbol{g},t\partial g({\boldsymbol{x}}))&=\operatorname{dist}^2\left(\boldsymbol{g},t\lambda_1\partial\|\cdot\|_1({\boldsymbol{x}}) +t\lambda_2\boldsymbol{D}^{\top}\partial\|\cdot\|_1(\boldsymbol{Dx})\right)
		\leq\|\boldsymbol{g}-t(\lambda_1\boldsymbol{z}_1+\lambda_2\boldsymbol{D}^{\top}\boldsymbol{z}_2)\|_2^2\nonumber\\
		&=\|\boldsymbol{g}-t\lambda_1(\operatorname{sign}({\boldsymbol{x}})_{\mathbb{S}_R}+\operatorname{sign}(\boldsymbol{g})_{\mathbb{S}_R^c})-t\lambda_2(\boldsymbol{D}^{\top}\operatorname{sign}(\boldsymbol{Dx})_{\mathbb{S}_G}+\boldsymbol{D}^{\top}\operatorname{sign}(\boldsymbol{Dg})_{\mathbb{S}_G^c})\|_2^2\nonumber\\
		&=\textup{\uppercase\expandafter{\romannumeral1}}+\textup{\uppercase\expandafter{\romannumeral2}}+\textup{\uppercase\expandafter{\romannumeral3}}\label{l1tv-36},
	\end{align}
	with	\begin{align*}
		\textup{\uppercase\expandafter{\romannumeral1}}=&\|\boldsymbol{g}\|_2^2+t^2\lambda_1^2(\|\operatorname{sign}({\boldsymbol{x}})_{\mathbb{S}_R}\|_2^2+\|\operatorname{sign}(\boldsymbol{g})_{\mathbb{S}_R^c}\|_2^2)+t^2\lambda_2^2(\|\boldsymbol{D}^{\top}\operatorname{sign}(\boldsymbol{Dx})_{\mathbb{S}_G}\|_2^2+\|\boldsymbol{D}^{\top}\operatorname{sign}(\boldsymbol{Dg})_{\mathbb{S}_G^c}\|_2^2),\\
		\textup{\uppercase\expandafter{\romannumeral2}}=&2t^2\lambda_1\lambda_2(\langle\operatorname{sign}(\boldsymbol{g})_{\mathbb{S}_R^c},\boldsymbol{D}^{\top}\operatorname{sign}(\boldsymbol{Dg})_{\mathbb{S}_G^c}\rangle+\langle\operatorname{sign}({\boldsymbol{x}})_{\mathbb{S}_R},\boldsymbol{D}^{\top}\operatorname{sign}(\boldsymbol{Dx})_{\mathbb{S}_G}\rangle) \\ &-2t(\lambda_1\langle\boldsymbol{g},\operatorname{sign}(\boldsymbol{g})_{\mathbb{S}_R^c}\rangle+\lambda_2\langle\boldsymbol{g},\boldsymbol{D}^{\top}\operatorname{sign}(\boldsymbol{Dg})_{\mathbb{S}_G^c}\rangle),\\
		\textup{\uppercase\expandafter{\romannumeral3}}=&2t^2\lambda_1\lambda_2(\langle\operatorname{sign}({\boldsymbol{x}})_{\mathbb{S}_R},\boldsymbol{D}^{\top}\operatorname{sign}(\boldsymbol{Dg})_{\mathbb{S}_G^c}\rangle+\langle\operatorname{sign}(\boldsymbol{g})_{\mathbb{S}_R^c},\boldsymbol{D}^{\top}\operatorname{sign}(\boldsymbol{Dx})_{\mathbb{S}_G}\rangle) \\
		&-2t^2(\lambda_1^2\langle\operatorname{sign}({\boldsymbol{x}})_{\mathbb{S}_R},\operatorname{sign}(\boldsymbol{g})_{\mathbb{S}_R^c}\rangle+\lambda_2^2\langle \boldsymbol{D}^{\top}\operatorname{sign}(\boldsymbol{Dx})_{\mathbb{S}_G},\boldsymbol{D}^{\top}\operatorname{sign}(\boldsymbol{Dg})_{\mathbb{S}_G^c}\rangle) \\
		&-2t(\lambda_1\langle\boldsymbol{g},\operatorname{sign}({\boldsymbol{x}})_{\mathbb{S}_R}\rangle+\lambda_2\langle\boldsymbol{g},\boldsymbol{D}^{\top}\operatorname{sign}(\boldsymbol{Dx})_{\mathbb{S}_G}\rangle).
	\end{align*}
	Taking expectation on both sides of \eqref{l1tv-36} gives
	\begin{equation}\label{l1tv-38}
		\operatorname{E}[\operatorname{dist}^2(\boldsymbol{g},t\partial g({\boldsymbol{x}}))]\leq\operatorname{E}[\textup{\uppercase\expandafter{\romannumeral1}}]+
		\operatorname{E}[\textup{\uppercase\expandafter{\romannumeral2}}]+\operatorname{E}[\textup{\uppercase\expandafter{\romannumeral3}}].
	\end{equation}
	Next, we bound the terms in
	\eqref{l1tv-38}.
	For the term $\operatorname{E}[\textup{\uppercase\expandafter{\romannumeral1}}]$, we have
	\begin{align}\label{l1tv-n5}
		\operatorname{E}[{\rm I}]\overset{(a)}{=}&n+t^2\Bigg(\lambda_1^2(|\mathbb{S}_R|+|\mathbb{S}_R^c|)+\lambda_2^2\sum_{j,k\in \mathbb{S}_G}\boldsymbol{d}_j^{\top}\boldsymbol{d}_k\operatorname{sign}(\boldsymbol{d}_j^{\top}{\boldsymbol{x}}) \operatorname{sign}(\boldsymbol{d}_k^{\top}{\boldsymbol{x}})\nonumber\\
		&+\lambda_2^2\sum\limits_{j,k\in \mathbb{S}_G^c}\boldsymbol{d}_j^{\top}\boldsymbol{d}_k\frac{2}{\pi}\arcsin\frac{\boldsymbol{d}_j^{\top}\boldsymbol{d}_k }{\|\boldsymbol{d}_j\|_2\|\boldsymbol{d}_k\|_2}\Bigg)\nonumber\\
		\overset{(b)}{\leq}& n+t^2\left(\lambda_1^2(|\mathbb{S}_R|+|\mathbb{S}_R^c|)+\lambda_2^2\left(2|\mathbb{S}_G|+2s_1+2|\mathbb{S}_G^c|+\frac{2}{3}s_2\right)\right)\nonumber\\
		\leq& n+t^2\left(\lambda_1^2(|\mathbb{S}_R|+|\mathbb{S}_R^c|)+\lambda_2^2\left(4|\mathbb{S}_G|+\frac{8}{3}|\mathbb{S}_G^c|-\frac{8}{3}\right)\right),
	\end{align}
	where equality ($a$) holds, due to the following identities
	\begin{equation}\label{l1tv-n2}
		\begin{cases}
			\operatorname{E}[\|\boldsymbol{g}\|_2^2]=n,\quad
			\operatorname{E}[\|\operatorname{sign}({\boldsymbol{x}})_{\mathbb{S}_R}\|_2^2]=|\mathbb{S}_R|, \quad
			\operatorname{E}[\|\operatorname{sign}(\boldsymbol{g})_{\mathbb{S}_R^c}\|_2^2]=|\mathbb{S}_R^c|,\\
			\operatorname{E}[\|\boldsymbol{D}^{\top}\operatorname{sign}(\boldsymbol{Dx})_{\mathbb{S}_G}\|_2^2]=\sum\limits_{j,k\in \mathbb{S}_G}\boldsymbol{d}_j^{\top}\boldsymbol{d}_k\operatorname{sign}(\boldsymbol{d}_j^{\top}{\boldsymbol{x}})\operatorname{sign}(\boldsymbol{d}_k^{\top}{\boldsymbol{x}}),\\
			\operatorname{E}[\|\boldsymbol{D}^{\top}\operatorname{sign}(\boldsymbol{Dg})_{\mathbb{S}_G^c}\|_2^2]=\sum\limits_{j,k\in \mathbb{S}_G^c}\boldsymbol{d}_j^{\top}\boldsymbol{d}_k\operatorname{sign}(\boldsymbol{d}_j^{\top}\boldsymbol{g})\operatorname{sign}(\boldsymbol{d}_k^{\top}\boldsymbol{g})\overset{(c)}{=}\sum\limits_{j,k\in \mathbb{S}_G^c}\boldsymbol{d}_j^{\top}\boldsymbol{d}_k\frac{2}{\pi}\arcsin\frac{\boldsymbol{d}_j^{\top}\boldsymbol{d}_k}{\|\boldsymbol{d}_j\|_2\|\boldsymbol{d}_k\|_2},
		\end{cases}
	\end{equation}
	where (c) is due to Lemma \ref{l1tv-lemma2}. The inequality ($b$) holds since $\boldsymbol{d}_j^{\top}\boldsymbol{d}_k=-1$ if $|j-k|=1$, $\|\boldsymbol{d}_j\|_2=\sqrt{2}$ if $j\in[n-1]$ and $\boldsymbol{d}_j^{\top}\boldsymbol{d}_k=0$ otherwise.	For the term $\operatorname{E}[\rm II]$, we have
	\begin{align}\label{l1tv-n3}
		\operatorname{E}[{\rm II}]=&2t^2\lambda_1\lambda_2\Bigg(\sum\limits_{i\in \mathbb{S}_G^c}\operatorname{E}\left[\boldsymbol{d}_i^{\top}\operatorname{sign}(\boldsymbol{g})_{\mathbb{S}_R^c}\times\operatorname{sign}(\boldsymbol{d}_i^{\top}\boldsymbol{g})\right]+\sum\limits_{i\in \mathbb{S}_G}\boldsymbol{d}_i^{\top}\operatorname{sign}({\boldsymbol{x}})_{\mathbb{S}_R}\times\operatorname{sign}(\boldsymbol{d}_i^{\top}{\boldsymbol{x}})\Bigg)\nonumber\\ &-2t\Bigg(\lambda_1\sqrt{\frac{2}{\pi}}|\mathbb{S}_R^c|+\lambda_2\sqrt{\frac{2}{\pi}}\sum\limits_{i\in \mathbb{S}_G^c}\|\boldsymbol{d}_i\|_2\Bigg)\nonumber\\
		\overset{(d)}{\leq}& 2t^2\lambda_1\lambda_2\big(0+2\min\{|\mathbb{S}_R|,|\mathbb{S}_G|\}\big)-2\sqrt{\frac{2}{\pi}}t\left(\lambda_1|\mathbb{S}_R^c|+\lambda_2\sqrt{2}|\mathbb{S}_G^c|\right)\nonumber\\
		\leq&4t^2\lambda_1\lambda_2\min\{|\mathbb{S}_R|,|\mathbb{S}_G|\}-2\sqrt{\frac{2}{\pi}}t\left(\lambda_1|\mathbb{S}_R^c|+\sqrt{2}\lambda_2|\mathbb{S}_G^c|\right),
	\end{align}
	where $(d)$ holds by Lemma \ref{l1tv-lemma9} with $\{i,i+1\}\subset\mathbb{S}_R^c$ or $\{i,i+1\}\subset\mathbb{S}_R$ on the condition $i\in\mathbb{S}_G^c$.
	For the term $\operatorname{E}[{\rm III}]$, we have
	\begin{equation}\label{l1tv-n4}
		\operatorname{E}[{\rm III}]=0.
	\end{equation}
	By substituting \eqref{l1tv-n5}, \eqref{l1tv-n3} and \eqref{l1tv-n4} into \eqref{l1tv-38},
	and noting $|\mathbb{S}_R|=s_r$, $|\mathbb{S}_G|=s_g$, $|\mathbb{S}_R^c|=s_r-1$, $|\mathbb{S}_G^c|=n-1-s_g$, we get
	\begin{align}\label{l1tv-39}
		\operatorname{dist}^2(\boldsymbol{g},t\partial g({\boldsymbol{x}}))
		\leq&\left(\lambda_1^2(|\mathbb{S}_R|+|\mathbb{S}_R^c|)+\lambda_2^2\left(4|\mathbb{S}_G|+\frac{8}{3}|\mathbb{S}_G^c|-\frac{8}{3}\right)+4\lambda_1\lambda_2\min\{|\mathbb{S}_R|,|\mathbb{S}_G|\}\right)t^2\nonumber\\
		&-2\sqrt{\frac{2}{\pi}}\left(\lambda_1|\mathbb{S}_R^c|+\sqrt{2}\lambda_2|\mathbb{S}_G^c|\right)t+n\nonumber\\
		\leq&\left(\lambda_1^2n+\lambda_2^2\left(\frac{4}{3}s_g+\frac{8}{3}n-\frac{16}{3}\right)+4\lambda_1\lambda_2\min\{s_r,s_g\}\right)t^2\nonumber\\
		&-2\sqrt{\frac{2}{\pi}}\left(\lambda_1(n-s_r)+\sqrt{2}\lambda_2(n-1-s_g)\right)t+n.
	\end{align}
	Meanwhile, by \cite{daei2019on}, we have
	\begin{equation}\label{l1tv-n15}
		\delta(\mathcal{D}(g,{\boldsymbol{x}})) \leq \inf_{t\geq0}\ \operatorname{E}[\ \textup{dist}^2(\boldsymbol{g},t\partial g({\boldsymbol{x}}))].
	\end{equation}
	By minimizing \eqref{l1tv-39} in $t$ and combining with \eqref{l1tv-n15}, we get
	\begin{equation*}
		\delta(\mathcal{D}(g,{\boldsymbol{x}}))\leq n-\frac{6}{\pi}\frac{\left[\lambda_1(n-s_r)+\sqrt{2}\lambda_2(n-1-s_g)\right]^2}{3n\lambda_1^2+4\left(2n+s_g-4\right)\lambda_2^2+12\lambda_1\lambda_2\min\{s_r,s_g\}}=\Phi(s_r,s_g).
	\end{equation*}
	This completes the proof of the theorem.
	
	\subsection{Proof of Theorem \ref{l1tv-theorem3.2}}\label{l1tv-proof2}
	
	\begin{proof}
		By \cite[Corollary 3.5]{tropp2015convex}, we have that with a probability at least $1-e^{-\frac{t^2}{2}}$, there holds
		\begin{equation}\label{l1tv-14}
			\|\boldsymbol{x}^\ast - \hat{\boldsymbol{x}}\|_2\leq \frac{2\epsilon}{[\sqrt{m-1}-w(\mathcal{D}(g,\boldsymbol{x}^\ast))-t]_+}.
		\end{equation}
		In addition, Remark \ref{l1tv-remark11} of Theorem \ref{l1tv-theorem3.1} implies  $w(\mathcal{D}(g,\boldsymbol{x}^\ast))\leq \sqrt{\Phi(s_r,s_g)}$.
		By	substituting it into \eqref{l1tv-14}, and combining with the assumption $m>(\sqrt{\Phi(s_r,s_g)}+t)^2+1$, we complete the proof.
	\end{proof}
	
	\subsection{Proof of Theorem \ref{l1tv-thm2}}\label{l1tv-proof3}
	First we outline the overall proof strategy. For part (a), we first show that there exists $(u,t)$ such that $\mathcal{G}_{1/u}^{f,g_1}({\boldsymbol{x}}^k(u,t)) \in \partial F_1({\boldsymbol{x}}^k(u,t))$, and
	then Algorithm \ref{alg:pgm-ista} can be viewed as the proximal subgradient method \cite{beck2017first},
	whose global convergence is guaranteed. Part (b) is proved by showing that the sequence
	$\{{\boldsymbol{x}}^k\}_{k\geq 0}$ generated by Algorithm \ref{alg:pgm-ista} is Fej\'er monotone,
	which converges to some fixed points. First we give a useful lemma.
	
	\begin{lemma}\label{l1tv-lemma6}
		Let $\boldsymbol{x}\in\mathbb{R}^n$, $\boldsymbol{A}=[\boldsymbol{a}_1,\boldsymbol{a}_2,\cdots,\boldsymbol{a}_n]$ and $\mathbb{I}_{\neq}(\boldsymbol{x}) = \{i|x_i\neq 0, i\in[n]\}$.
		Suppose that $u_0>0$ satisfies
		\begin{equation}\label{l1tv-16}
			\frac{1}{u_0}\geq\max\left\{\frac{\operatorname{sign}(x_i)\lambda_1+\sum_{j=1}^{n}\boldsymbol{a}_i^{\top}\boldsymbol{a}_jx_j-\boldsymbol{a}_i^{\top}\boldsymbol{y}}{x_i}, i\in \mathbb{I}_{\neq}(\boldsymbol{x})\right\}.
		\end{equation}
		Then, for any $u\in(0,u_0]$,
		\begin{equation}\label{l1tv-17}
			\mathcal{G}_{1/u}^{f,g_1}(\boldsymbol{x}) \subset \partial F_1(\boldsymbol{x}),
		\end{equation}
		where $F_1(\boldsymbol{x})=f(\boldsymbol{x})+g_1(\boldsymbol{x})$ with $f(\boldsymbol{x})=\frac{1}{2}\|\boldsymbol{y}-A\boldsymbol{x}\|_2^2$ and $g_1(\boldsymbol{x})=\lambda_1\|\boldsymbol{x}\|_1$.
	\end{lemma}
	\begin{proof}
		The $i$-th coordinate of $\partial F_1({\boldsymbol{x}})$ is
		\begin{equation}\label{l1tv-45}
			\left[\partial F_1({\boldsymbol{x}})\right]_i\in\begin{cases}
				[\nabla f({\boldsymbol{x}})]_i + \lambda_1, &{x_i>0}\\
				[\nabla f({\boldsymbol{x}})]_i + [-\lambda_1,\lambda_1], &{x_i=0}\\
				[\nabla f({\boldsymbol{x}})]_i - \lambda_1, &{x_i<0}
			\end{cases}.
		\end{equation}
		Moreover,
		\begin{align*}
			\mathcal{G}_{1/u}^{f,g_1}({\boldsymbol{x}}) =& \tfrac{1}{u}\left({\boldsymbol{x}}-\operatorname{prox}_{ug_1}\left({\boldsymbol{x}}-u\nabla f({\boldsymbol{x}})\right)\right)\\
			=& \tfrac{1}{u}\left({\boldsymbol{x}}-\operatorname{sign}\left({\boldsymbol{x}}-u\nabla f({\boldsymbol{x}})\right)\odot\left[|{\boldsymbol{x}}-u\nabla f({\boldsymbol{x}})|-\lambda_1u\right]_+\right).
		\end{align*}
		Hence, the $i$-th coordinate of $\mathcal{G}_{1/u}^{f,g_1}({\boldsymbol{x}})$ is
		\begin{equation}\label{l1tv-46}
			\left[\mathcal{G}_{1/u}^{f,g_1}({\boldsymbol{x}})\right]_i\in\begin{cases}
				[\nabla f({\boldsymbol{x}})]_i + \lambda_1, &{\left[\frac{1}{u}{\boldsymbol{x}}-\nabla f({\boldsymbol{x}})\right]_i\geq\lambda_1}\\
				[\nabla f({\boldsymbol{x}})]_i + (-\lambda_1,\lambda_1), &{\left[\frac{1}{u}{\boldsymbol{x}}-\nabla f({\boldsymbol{x}})\right]_i\in(-\lambda_1,\lambda_1)}\\
				[\nabla f({\boldsymbol{x}})]_i - \lambda_1, &{\left[\frac{1}{u}{\boldsymbol{x}}-\nabla f({\boldsymbol{x}})\right]_i\leq\lambda_1}
			\end{cases}.
		\end{equation}
		Thus, if $x_i \in \mathbb{I}_0({\boldsymbol{x}}) = \{i|x_i= 0, i\in[n]\}$, it follows from \eqref{l1tv-45} and \eqref{l1tv-46} that
		\begin{equation}\label{l1tv-47}
			\left[\mathcal{G}_{1/u}^{f,g_1}({\boldsymbol{x}})\right]_i\subset\left[\partial F_1({\boldsymbol{x}})\right]_i.
		\end{equation}
		If $x_i > 0$, it follows from \eqref{l1tv-16} that
		\begin{equation*}
			\frac{1}{u}\geq\frac{1}{u_0}\geq\frac{\lambda_1+\sum_{j=1}^{n}\boldsymbol{a}_i^{\top}\boldsymbol{a}_jx_j-\boldsymbol{a}_i^{\top}\boldsymbol{y}}{x_i},
		\end{equation*}
		which implies
		\begin{equation}\label{l1tv-48}
			\left[\tfrac{1}{u}{\boldsymbol{x}}-\nabla f({\boldsymbol{x}})\right]_i\geq\lambda_1.
		\end{equation}
		Combining \eqref{l1tv-45} and \eqref{l1tv-48} gives
		\begin{equation}\label{l1tv-49}
			\left[\mathcal{G}_{1/u}^{f,g_1}({\boldsymbol{x}})\right]_i=\left[\partial F_1({\boldsymbol{x}})\right]_i.
		\end{equation}	
		Similarly, if $x_i < 0$, repeating the argument for $x_i>0$ implies that \eqref{l1tv-49} also holds for $u\in(0,u_0]$.
		Combining \eqref{l1tv-47} and \eqref{l1tv-49} gives \eqref{l1tv-17}.	
	\end{proof}
	
	Lemma \ref{l1tv-lemma6} indicates $\mathcal{G}_{1/u}^{f,g_1}(\boldsymbol{x}) \subset \partial F_1(\boldsymbol{x})$ for suitable $u$, connecting the proposed algorithms with the proximal subgradient method \cite{beck2017first}.
	
	\textbf{Proof of Theorem \ref{l1tv-thm2}.}
	\textbf{Part (a)}. We first claim that for any integer $k\geq0$, there exist $u>0$ and $t>0$, such that the sequence $\{{\boldsymbol{x}}^k(u,t)\}_{k\geq 0}$ satisfies
	\begin{equation}\label{l1tv-50}
		\mathcal{G}_{1/u}^{f,g_1}({\boldsymbol{x}}^k(u,t)) \in \partial F_1({\boldsymbol{x}}^k(u,t)).
	\end{equation}
	The proof of \eqref{l1tv-50} will be given later. The relation \eqref{l1tv-50} implies that Algorithm \ref{alg:pgm-ista} can be viewed
	as the proximal subgradient method \cite{beck2017first}.
	Following the theory of the proximal subgradient method, there exists a constant step size
	$t>0$ satisfying ${\boldsymbol{x}}(u,t) = {\boldsymbol{x}}^\ast$ and
	\begin{equation}\label{l1tv-51}
		{\boldsymbol{x}}(u,t) = \operatorname{prox}_{tg_2}\big({\boldsymbol{x}}(u,t)-t\mathcal{G}_{1/u}^{f,g_1}({\boldsymbol{x}}(u,t))\big).
	\end{equation}
	Next, we prove \eqref{l1tv-50} by mathematical induction. For $k=0,1$, by Lemma \ref{l1tv-lemma6}, there exist $u_0, u_1\in(0,\infty)$ satisfying \eqref{l1tv-16} for ${\boldsymbol{x}}^0$ and ${\boldsymbol{x}}^1$ respectively, and for any $u'\in(0,u_0]$, $u''\in(0,u_1]$,
	\begin{equation*}
		\mathcal{G}_{1/u'}^{f,g_1}({\boldsymbol{x}}^0)\in\partial F_1({\boldsymbol{x}}^0)\quad \textup{and} \quad \mathcal{G}_{1/u''}^{f,g_1}({\boldsymbol{x}}^1)\in\partial F_1({\boldsymbol{x}}^1).
	\end{equation*}
	Choosing $u=\min\{u_0,u_1\}$ leads to
	\begin{equation*}
		\mathcal{G}_{1/u}^{f,g_1}({\boldsymbol{x}}^0)\in\partial F_1({\boldsymbol{x}}^0)\quad \textup{and} \quad \mathcal{G}_{1/u}^{f,g_1}({\boldsymbol{x}}^1)\in\partial F_1({\boldsymbol{x}}^1).
	\end{equation*}
	Suppose that for an integer $k\geq0$, $\mathcal{G}_{1/u}^{f,g_1}(\{{\boldsymbol{x}}^k)\in\partial F_1(\{{\boldsymbol{x}}^k)$, by induction, the above arguments indicate that we can always find $u>0$ so that
	$\mathcal{G}_{1/u}^{f,g_1}(({\boldsymbol{x}}^{k+1})\in\partial F_1({\boldsymbol{x}}^{k+1})$,
	which implies the claim \eqref{l1tv-50}.\\
	
	\noindent\textbf{Part (b)}. The proof of part (a) indicates that a small parameter pair $(u,t)$ can ensure the convergence.
	Hence, we can choose $(u,t)$ with $u\in(0,2/\|\boldsymbol{A}\|_2^2)$ and $t\in(0,u]$ which also implies the convergence of Algorithm \ref{alg:pgm-ista}.
	That is, \eqref{l1tv-51} holds for some $\boldsymbol{x}(u_0,t_0)={\boldsymbol{x}}^\ast$ with $u_0\in(0,2/\|\boldsymbol{A}\|_2^2)$ and $t_0\in(0,u_0]$.
	For any $k \geq 0$, upon letting ${\boldsymbol{x}}^k(u,t)={\boldsymbol{x}}^k$, we have
	\begin{align}
		\label{l1tv-fixpoint}
		\left\|{\boldsymbol{x}}^{k+1}-{\boldsymbol{x}}^\ast\right\|_2=&\left\|\operatorname{prox}_{tg_2}\left({\boldsymbol{x}}^{k}-t\mathcal{G}_{1/u}^{f,g_1}\left({\boldsymbol{x}}^{k}\right)\right)-\operatorname{prox}_{tg_2}\left({\boldsymbol{x}}^\ast-t\mathcal{G}_{1/u}^{f,g_1}\left({\boldsymbol{x}}^\ast\right)\right)\right\|_2\nonumber\\
		\stackrel{(a)}{\leq}&\left\|\left({\boldsymbol{x}}^{k}-{\boldsymbol{x}}^\ast\right)-t\left(\mathcal{G}_{1/u}^{f,g_1}\left({\boldsymbol{x}}^{k}\right)-\mathcal{G}_{1/u}^{f,g_1}\left({\boldsymbol{x}}^\ast\right)\right)\right\|_2\nonumber\\
		\stackrel{(b)}{\leq}&\frac{t}{u}\left\|\operatorname{prox}_{ug_1}\left({\boldsymbol{x}}^{k}-u\boldsymbol{A}^{\top}(\boldsymbol{A}{\boldsymbol{x}}^{k}-\boldsymbol{y})\right)-\operatorname{prox}_{ug_1}\left({\boldsymbol{x}}^\ast-u\boldsymbol{A}^{\top}(\boldsymbol{A}{\boldsymbol{x}}^\ast-\boldsymbol{y})\right)\right\|_2\nonumber\\
		&+\left|1-\frac{t}{u}\right|\cdot\left\|{\boldsymbol{x}}^{k}-{\boldsymbol{x}}^\ast\right\|_2\nonumber\\
		\stackrel{(c)}{\leq}&\frac{t}{u}\left\|(\boldsymbol{I}_n-u\boldsymbol{A}^{\top}\boldsymbol{A})({\boldsymbol{x}}^{k}-{\boldsymbol{x}}^\ast)\right\|_2^2+\left|1-\frac{t}{u}\right|\cdot\left\|{\boldsymbol{x}}^{k}-{\boldsymbol{x}}^\ast\right\|_2\nonumber\\
		\stackrel{(d)}{\leq}&\rho(u,t)\left\|{\boldsymbol{x}}^{k}-{\boldsymbol{x}}^\ast\right\|_2,
	\end{align}
	where $\rho(u,t)=|1-\frac{t}{u}|+\frac{t}{u}\|\boldsymbol{I}_n-u\boldsymbol{A}^{\top}\boldsymbol{A}\|_2$, $(a)$ and $(c)$ hold by the nonexpansivity of the proximal operator \cite[Theorem 6.42]{beck2017first}, and $(b)$ and $(d)$ hold by the triangular inequality.
	Since $u\in(0,2/\|\boldsymbol{A}\|_2^2)$, $t\in(0,u]$, we get $\rho(u,t)=1-\frac{t}{u}+\frac{t}{u}\|\boldsymbol{I}_n-u\boldsymbol{A}^{\top}\boldsymbol{A}\|_2$.
	Further, since $\|\boldsymbol{I}_n-u\boldsymbol{A}^{\top}\boldsymbol{A}\|_2=\|u\boldsymbol{A}^{\top}\boldsymbol{A}-\boldsymbol{I}_n\|_2$, we obtain
	\begin{align*}
		&\|\boldsymbol{I}_n-u\boldsymbol{A}^{\top}\boldsymbol{A}\|_2=\max_{{\boldsymbol{x}}\in \mathbb{S}^{n-1}}\langle(\boldsymbol{I}_n-u\boldsymbol{A}^{\top}\boldsymbol{A}){\boldsymbol{x}},{\boldsymbol{x}}\rangle\\
		=&\max_{{\boldsymbol{x}}\in \mathbb{S}^{n-1}}\left|1-u\|\boldsymbol{A}{\boldsymbol{x}}\|_2^2\right|
		=\max\{|1-us_n^2|,|1-us_1^2|\},
	\end{align*}
	where $s_1$ and $s_n$ are the maximal and minimal singular values of $\boldsymbol{A}$ and the last equality holds since $s_n^2\leq\|\boldsymbol{Ax}\|_2^2\leq s_1^2$.
	Since $\boldsymbol{A}\in\mathbb{R}^{m\times n}$ with $m<n$, we obtain $s_n=0$.
	Thus for $u\in(0,2/\|\boldsymbol{A}\|_2^2)$, $\|\boldsymbol{I}_n-u\boldsymbol{A}^{\top}\boldsymbol{A}\|_2=\max\{1,|1-us_1^2|\} = 1$.
	These discussions imply $\rho(u,t)=1$.
	By substituting the identity into \eqref{l1tv-fixpoint}, we get
	\begin{equation}\label{l1tv-54}
		\left\|{\boldsymbol{x}}^{k+1}-{\boldsymbol{x}}^\ast\right\|_2\leq\left\|{\boldsymbol{x}}^{k}-{\boldsymbol{x}}^\ast\right\|_2.
	\end{equation}
	This estimate implies that the sequence $\{{\boldsymbol{x}}^k\}_{k\geq 0}$ generated by Algorithm \ref{alg:pgm-ista}
	is Fej\'er monotone and hence converges to some fixed points \cite[Theorem 8.16]{beck2017first}.
	\qed
	
	\subsection{Proof of Theorem \ref{l1tv-thm4}}\label{l1tv-proof4}
	First we describe the overall proof strategy. We first construct a $4/(3u)$-smooth function $H_u(\cdot)$ to obtain $F_u(\cdot)=H_u(\cdot)+g_2(\cdot)$.
	Then $F_u(\cdot)$ is used to approximate $F(\cdot)$, which yields an upper bound of $F({\boldsymbol{x}}^{k+1}) - F(\boldsymbol{x}^\ast)$.
	Below we denote $f(\boldsymbol{x}) = \frac{1}{2}\boldsymbol{x}^{\top}\boldsymbol{Q}\boldsymbol{x} + \boldsymbol{b}^{\top}\boldsymbol{x} +c$, with $\boldsymbol{Q} = \boldsymbol{A}^{\top}\boldsymbol{A}$, $\boldsymbol{b} = \boldsymbol{A}^{\top}\boldsymbol{y}$ and $c = \frac{1}{2} \boldsymbol{y}^{\top}\boldsymbol{y}$.
	We denote by $l_f$, $l_{g_1}$ and $l_{g_2}$ the Lipschitz constants of $f(\cdot)$, $g_1(\cdot)$ and $g_2(\cdot)$.
	To prove Theorem \ref{l1tv-thm4}, we need the next lemma on the gradient mapping operator.
	
	\begin{lemma}[{\cite[Lemma 2.1]{sulam2020on}}]\label{l1tv-lemma8}
		For any $u>0$ and ${\boldsymbol{x}}\in\mathbb{R}^n$,
		\begin{equation*}
			\left\|\mathcal{G}_{1/u}^{f,g_1}({\boldsymbol{x}})\right\|_2\leq \|\boldsymbol{Q}\|_2r + \|\boldsymbol{b}\|_2 + l_{g_1}.
		\end{equation*}
	\end{lemma}
	
	Now, we are ready to prove Theorem \ref{l1tv-thm4}.
	
	\textbf{Proof of Theorem \ref{l1tv-thm4}.}
	Consider the function $H_u :\mathbb{R}^n \rightarrow \mathbb{R}$ given by
	\begin{equation*}
		H_u({\boldsymbol{x}}) = \tfrac{1}{2}{\boldsymbol{x}}^{\top}(\boldsymbol{Q}-u\boldsymbol{Q}^2){\boldsymbol{x}} - \boldsymbol{b}^{\top}(\boldsymbol{I}_n-u\boldsymbol{Q}){\boldsymbol{x}} + M_{g_1}^u\left((\boldsymbol{I}_n-u\boldsymbol{Q}){\boldsymbol{x}} + u\boldsymbol{b}\right) + \tfrac{1}{2} \boldsymbol{y}^{\top}\boldsymbol{y},
	\end{equation*}
	where $M_{g_1}^u$ is the Moreau envelope of $g_1$ with a smoothness parameter $u$ (see \cite{moreau1965proximite} for Moreau envelope).
	Note that $H_u$ is convex since $u\in(0,2/\|\boldsymbol{Q}\|_2)$ and the Moreau envelope of a convex function is convex.
	Since the gradient of the Moreau envelop is given by $\nabla M_{g_1}^u({\boldsymbol{x}}) = \frac{1}{u}({\boldsymbol{x}}-\operatorname{prox}_{ug_1}({\boldsymbol{x}}))$, we have
	\begin{equation}\label{l1tv-55}
		\nabla H_u({\boldsymbol{x}}) = (\boldsymbol{I}_n-u\boldsymbol{Q})\mathcal{G}_{1/u}^{f,g_1}({\boldsymbol{x}}).
	\end{equation}
	Let $\boldsymbol{a}_1 = \frac{1}{t}({\boldsymbol{x}}^k - \operatorname{prox}_{tg_2}({\boldsymbol{x}}^k - t \nabla H_u({\boldsymbol{x}}^k))) = \mathcal{G}_{1/t}^{H_u,g_2}({\boldsymbol{x}}^k)$ and
	$\boldsymbol{a}_2 = \frac{1}{t}({\boldsymbol{x}}^k - \operatorname{prox}_{tg_2}({\boldsymbol{x}}^k - t \mathcal{G}_{1/u}^{f,g_1}({\boldsymbol{x}}^k))) = \frac{1}{t}({\boldsymbol{x}}^k - {\boldsymbol{x}}^{k+1})$.
	Then by the triangle inequality, we have
	\begin{equation}\label{l1tv-56}
		\|\boldsymbol{a}_1\|_2\leq\|\boldsymbol{a}_2\|_2 + \|\boldsymbol{a}_1 - \boldsymbol{a}_2\|_2.
	\end{equation}
	Next, by \eqref{l1tv-55} and nonexpansivity of the proximal operator  \cite[Theorem 6.42]{beck2017first}, and Lemma \ref{l1tv-lemma8}, we have
	\begin{align}\label{l1tv-57}
		\|\boldsymbol{a}_1 - \boldsymbol{a}_2\|_2 =& \frac{1}{t}\left\| \operatorname{prox}_{tg_2}\left({\boldsymbol{x}}^k - t \nabla H_u({\boldsymbol{x}}^k)\right)-\operatorname{prox}_{tg_2}\left({\boldsymbol{x}}^k - t \mathcal{G}_{1/u}^{f,g_1}({\boldsymbol{x}}^k)\right)\right\|_2\nonumber\\
		{\leq}& \left\|(\boldsymbol{I}_n-u\boldsymbol{Q})\mathcal{G}_{1/u}^{f,g_1}({\boldsymbol{x}}^k) - \mathcal{G}_{1/u}^{f,g_1}({\boldsymbol{x}}^k)\right\|_2
		= u\left\|\boldsymbol{Q}\mathcal{G}_{1/u}^{f,g_1}({\boldsymbol{x}}^k)\right\|_2\nonumber\\
		{\leq}& u\|\boldsymbol{Q}\|_2(\|\boldsymbol{Q}\|_2r + \|\boldsymbol{b}\|_2+l_{g_1})?
	\end{align}
	By assumption \eqref{l1tv-85}, we have
	\begin{equation}\label{l1tv-58}
		\|\boldsymbol{a}_2\|_2 =\tfrac{1}{t} \|{\boldsymbol{x}}^k - {\boldsymbol{x}}^{k+1}\|_2\leq \epsilon.
	\end{equation}
	By combining \eqref{l1tv-56}, \eqref{l1tv-57} and \eqref{l1tv-58}, we get
	\begin{equation}\label{l1tv-59}
		\|\boldsymbol{a}_1\|_2 = \|\mathcal{G}_{1/t}^{H_u,g_2}({\boldsymbol{x}}^k)\|_2\leq\epsilon + u\|\boldsymbol{Q}\|_2(\|\boldsymbol{Q}\|_2r + \|\boldsymbol{b}\|_2+l_{g_1}).
	\end{equation}
	Further, for any ${\boldsymbol{x}}_1$, ${\boldsymbol{x}}_2\in\mathbb{R}^n$, we have
	\begin{align*}
		\|\nabla H_u({\boldsymbol{x}}_1) - \nabla H_u({\boldsymbol{x}}_2)\|_2 =& \left\|(\boldsymbol{I}_n-u\boldsymbol{Q})\left(\mathcal{G}_{1/u}^{f,g_1}({\boldsymbol{x}}_1) - \mathcal{G}_{1/u}^{f,g_1}({\boldsymbol{x}}_2)\right)\right\|_2\nonumber\\
		\leq&\left\|\mathcal{G}_{1/u}^{f,g_1}({\boldsymbol{x}}_1) - \mathcal{G}_{1/u}^{f,g_1}({\boldsymbol{x}}_2)\right\|_2
		\leq\frac{4}{3u}\|{\boldsymbol{x}}_1-{\boldsymbol{x}}_2\|_2,
	\end{align*}
	where the first line holds since $u\in(0,2/\|\boldsymbol{Q}\|_2)$ and the second line holds since gradient mapping operators are firmly non-expansive with constant $3u/4$ \cite[Lemma 10.11]{beck2017first} (i.e., $H_u(\cdot)$ is $4/(3u)$-smooth, cf. Definition \ref{l1tv-lsmooth} for $l$-smooth).	Denote $F_u({\boldsymbol{x}}) = H_u({\boldsymbol{x}}) + g_2({\boldsymbol{x}})$ and a minimizer by ${\boldsymbol{x}}_u^\ast\in\mathop{\arg\min}_{\boldsymbol{x}} F_u({\boldsymbol{x}})$.
	Moreover, define
	\begin{equation*}
		\hat{{\boldsymbol{x}}}_u = \operatorname{prox}_{tg_2}\left({\boldsymbol{x}}^k - t \nabla H_u({\boldsymbol{x}}^k)\right).
	\end{equation*}
	By the fundamental prox-grad inequality \cite[Theorem 10.16]{beck2017first}, the convexity of $H_u(\cdot)$ and the three-points lemma  \cite{chen1993convergence}, since $t\in(0,3u/4)$, it follows that
	\begin{align*}
		F_u({\boldsymbol{x}}_u^\ast) - F_u(\hat{{\boldsymbol{x}}}_u) \geq& \frac{1}{2t}\|{\boldsymbol{x}}_u^\ast-\hat{{\boldsymbol{x}}}_u\|_2^2 - \frac{1}{2t}\|{\boldsymbol{x}}_u^\ast-{\boldsymbol{x}}^k\|_2^2\nonumber\\
		=&\frac{1}{2t}\left(2\langle {\boldsymbol{x}}^k-\hat{{\boldsymbol{x}}}_u,{\boldsymbol{x}}_u^\ast-\hat{{\boldsymbol{x}}}_u\rangle-\|{\boldsymbol{x}}^k-\hat{{\boldsymbol{x}}}_u\|_2^2\right),
	\end{align*}
	Thus, by the Cauchy-Schwarz inequality and the estimate $\|{\boldsymbol{x}}^{k}-{\boldsymbol{x}}_u^\ast\|_2\leq 2r$,
	\begin{align}\label{l1tv-61}
		F_u(\hat{{\boldsymbol{x}}}_u) - F_u({\boldsymbol{x}}_u^\ast) \leq& \frac{1}{2t}\|{\boldsymbol{x}}^k-\hat{{\boldsymbol{x}}}_u\|_2^2 + \frac{1}{t}\langle\hat{{\boldsymbol{x}}}_u- {\boldsymbol{x}}_u^\ast,{\boldsymbol{x}}^k-\hat{{\boldsymbol{x}}}_u\rangle\nonumber\\
		=&-\frac{1}{2t}\|{\boldsymbol{x}}^k-\hat{{\boldsymbol{x}}}_u\|_2^2 + \frac{1}{t}\langle {\boldsymbol{x}}^k-{\boldsymbol{x}}_u^\ast,{\boldsymbol{x}}^k-\hat{{\boldsymbol{x}}}_u\rangle\nonumber\\
		=&-\frac{1}{2t}\|{\boldsymbol{x}}^k-{\boldsymbol{x}}_u^\ast\|_2^2 + \langle {\boldsymbol{x}}^{k}-{\boldsymbol{x}}_u^\ast,\mathcal{G}_{1/t}^{H_u,g_2}({\boldsymbol{x}}^k)\rangle\nonumber\\
		\leq&\langle {\boldsymbol{x}}^{k}-{\boldsymbol{x}}_u^\ast,\mathcal{G}_{1/t}^{H_u,g_2}({\boldsymbol{x}}^k)\rangle
		\leq2r\left\|\mathcal{G}_{1/t}^{H_u,g_2}({\boldsymbol{x}}^k)\right\|_2
	\end{align}
	Combining \eqref{l1tv-59} with \eqref{l1tv-61} yields
	\begin{equation}\label{l1tv-62}
		F_u(\hat{{\boldsymbol{x}}}_u) - F_u({\boldsymbol{x}}_u^\ast) \leq 2r\epsilon + 2ur\|\boldsymbol{Q}\|_2(\|\boldsymbol{Q}\|_2r + \|\boldsymbol{b}\|_2+l_{g_1}).
	\end{equation}
	Finally, we connect $F_u({\boldsymbol{x}}_u^\ast)$, $F_u({\boldsymbol{x}}^k)$ with $F({\boldsymbol{x}}_u^\ast)$, $F({\boldsymbol{x}}^k)$, respectively. Note that for any ${\boldsymbol{x}}\in\mathbb{R}^n$, by basic properties of the Moreau envelope  \cite[Theorem 10.51]{beck2017first} in \eqref{l1tv-63} and the $l_{g_1}-$Lipschitz property of $g_1$, we have
	\begin{align}\label{l1tv-63}
		|H_u({\boldsymbol{x}}) - f({\boldsymbol{x}}) - g_1({\boldsymbol{x}})| =& \left|u\left(-\frac{1}{2}\|\boldsymbol{Q}{\boldsymbol{x}}\|_2^2+\boldsymbol{b}^{\top}\boldsymbol{Q}{\boldsymbol{x}}\right)+M_{g_1}^u\left((\boldsymbol{I}_n-u\boldsymbol{Q}){\boldsymbol{x}}+u\boldsymbol{b}\right)-g_1({\boldsymbol{x}})\right|\nonumber\\
		\leq& \left(\frac{1}{2}\|\boldsymbol{Q}\|_2^2r^2+\|\boldsymbol{b}\|_2\|\boldsymbol{Q}\|_2r\right)u+\left|M_{g_1}^u\left((\boldsymbol{I}_n-u\boldsymbol{Q}){\boldsymbol{x}}+u\boldsymbol{b}\right)-g_1({\boldsymbol{x}})\right|\nonumber\\
		\leq&\left(\frac{1}{2}\|\boldsymbol{Q}\|_2^2r^2+\|\boldsymbol{b}\|_2\|\boldsymbol{Q}\|_2r\right)u+\big|M_{g_1}^u\left((\boldsymbol{I}_n-u\boldsymbol{Q}){\boldsymbol{x}}+u\boldsymbol{b}\right)\nonumber\\
		&-g_1((\boldsymbol{I}_n-u\boldsymbol{Q}){\boldsymbol{x}}+u\boldsymbol{b})\big|+\left|g_1((\boldsymbol{I}_n-u\boldsymbol{Q}){\boldsymbol{x}}+u\boldsymbol{b})-g_1({\boldsymbol{x}})\right|\nonumber\\
		\leq&\left(\frac{1}{2}\|\boldsymbol{Q}\|_2^2r^2+\|\boldsymbol{b}\|_2\|\boldsymbol{Q}\|_2r\right)u+\frac{l_{g_1}^2}{2}u+l_{g_1}(\|\boldsymbol{Q}\|_2r + \|\boldsymbol{b}\|_2)u\\
		\leq&(\|\boldsymbol{Q}\|_2r+\|\boldsymbol{b}\|_2)(\|\boldsymbol{Q}\|_2r+l_{g_1})u+\frac{l_{g_1}^2}{2}u\nonumber.
	\end{align}
	Thus, we obtain
	$|H_u({\boldsymbol{x}}) - f({\boldsymbol{x}}) - g_1({\boldsymbol{x}})| = |F_u({\boldsymbol{x}})-F({\boldsymbol{x}})|\leq cu$,
	with $	c=(\|\boldsymbol{Q}\|_2r+\|\boldsymbol{b}\|_2)(\|\boldsymbol{Q}\|_2r+l_{g_1})+\frac{l_{g_1}^2}{2}$.
	Consequently,
	\begin{equation}
		\min_{{\boldsymbol{x}}} F({\boldsymbol{x}}) \geq \min_{{\boldsymbol{x}}}F_u({\boldsymbol{x}}) - cu,\quad \mbox{i.e.,}\quad \label{l1tv-64}
		F(\boldsymbol{x}^\ast) \geq F_u({\boldsymbol{x}}_u^\ast) - cu.
	\end{equation}
	Moreover, since $\|{\boldsymbol{x}}^{k+1}-\hat{{\boldsymbol{x}}}_u\|_2=t\|\boldsymbol{a}_1-\boldsymbol{a}_2\|_2$, we have	
	\begin{align*}
		F({\boldsymbol{x}}^{k+1}) - F(\hat{{\boldsymbol{x}}}_u) &\leq \left|f({\boldsymbol{x}}^{k+1}) - f(\hat{{\boldsymbol{x}}}_u)\right| + \left|g_1({\boldsymbol{x}}^{k+1}) - g_1(\hat{{\boldsymbol{x}}}_u)\right| + \left|g_2({\boldsymbol{x}}^{k+1}) - g_2(\hat{{\boldsymbol{x}}}_u)\right|\\
		&\leq\left(l_f+l_{g_1}+l_{g_2}\right)\|{\boldsymbol{x}}^{k+1}-\hat{{\boldsymbol{x}}}_u\|_2
		=\left(l_f+l_{g_1}+l_{g_2}\right)t\|\boldsymbol{a}_1-\boldsymbol{a}_2\|_2.
	\end{align*}
	Combining this with \eqref{l1tv-57} leads to
	\begin{equation}\label{l1tv-65}
		F({\boldsymbol{x}}^{k+1})\leq F(\hat{{\boldsymbol{x}}}_u) + t\kappa u,
	\end{equation}
	with $\kappa=\left(l_f+l_{g_1}+l_{g_2}\right)\|\boldsymbol{Q}\|_2(\|\boldsymbol{Q}\|_2r + \|\boldsymbol{b}\|_2+l_{g_1})$.
	Finally, combining \eqref{l1tv-62}, \eqref{l1tv-64} and \eqref{l1tv-65} gives
	\begin{align}\label{l1tv-87}
		F({\boldsymbol{x}}^{k+1}) - F(\boldsymbol{x}^\ast) &\leq 2r\epsilon + (2r\|\boldsymbol{Q}\|_2(\|\boldsymbol{Q}\|_2r + \|\boldsymbol{b}\|_2+l_{g_1}) + t\kappa +c)u =2r\epsilon + (t\kappa + \beta)u.
	\end{align}
	By \cite[Theorem 3.61]{beck2017first} (for Lipschitz functions), we can choose $l_f=\|\boldsymbol{Q}\|_2r+\|\boldsymbol{b}\|_2$.
	Obviously, $g_1(\cdot)$ and $g_2(\cdot)$ are Lipschitz with $l_{g_1} = \lambda_1$ and $l_{g_2} = \lambda_2\|\boldsymbol{D}\|_2$. Substituting them into \eqref{l1tv-87} completes the proof.
	\qed
	
	\section{Conclusion and discussions}\label{l1tv-section6}
	
	Motivated by the $\ell^1$ and TV penalties in compressed sensing, we have investigated the $\ell^1$-TV method, which simultaneously enforces sparsity and gradient sparsity in the signal. We have established the recovery guarantee of the $\ell^1$-TV method, which quantifies the required number of measurements depending on regular sparsity and gradient sparsity levels in order to achieve successful recovery.
	We have proposed the PGM-ISTA algorithm based on the gradient mapping and proximal operator, and analyzed theoretical properties of  PGM-ISTA, including global convergence and its parameter selections. In addition, we have proposed a learned solver for the regularized problem, by unrolling the PGM-ISTA. The experimental results on ECG signals show that LPGM-ISTA significantly outperforms traditional algorithms.
	
	Below we briefly discuss several topics that are pertinent to this work as well as future research topics.
	
	{\bf Combination of regularizers.}
	In many recovery problems, combined regularizers have shown enhanced recovery performance.
	The theoretical research, however, lags far behind, impeding its interpretability and practical applications.
	Our findings provide a robust recovery guarantee for the $\ell^1$-TV method, thereby stimulating further research on recovery problems using combined convex regularizers.
	
	{\bf Equivalence between constrained problem \eqref{l1tv-model} and the regularized problem \eqref{l1tv-un}.}
	Theoretical guarantees for the constrained form like \eqref{l1tv-model} seem preferred, but  numerical algorithms are mostly designed for the regularized counterparts like \eqref{l1tv-un}.
	Several recent studies have investigated the equivalence between the constrained problem and its regularized counterpart in some special cases \cite{li2022selecting,wang2021low}. However, existing studies fall far short in providing a general equivalence.
	Thus, it requires further research to expound the equivalence between problems \eqref{l1tv-model} and  \eqref{l1tv-un}.
	
	{\bf Relationship with fused LASSO.}
	The model \eqref{l1tv-un} resembles the fused LASSO (FLASSO) in statistics.
	Since combining $\ell^1$ and TV regularizers is efficient in inducing sparsity and gradient sparsity in signals, FLASSO has been widely applied in classification and coefficient selection \cite{li2014linearized,tibshirani2005sparsity,rapaport2008classification}.
	These successful applications confirm the benefit of combining $\ell^1$ and TV regularizers in the context of CS.
	The proposed LPGM-ISTA can be directly applied to solve the FLASSO problem.
	
	{\bf Theoretical perspectives of LPGM-ISTA.}
	In contrast to successful practical applications of algorithm unrolling, its theoretical understanding is still in its nascent stage.
	Recently, Chen et al. \cite{chen2018theoretical,liu2019alista} proved that LISTA can attain a linear convergence, which is better than the sublinear convergence of ISTA/FISTA in general scenarios. Their research findings can shed insights into the theoretical underpinnings of LPGM-ISTA.

	\section*{Acknowledgements}
	
	This research was funded by National Natural Science Foundation of China (Grants No. 12301594, 12071380), National Key Research and Development Program of China (Grant No. 2023YFA1008500), Joint Funds of the Natural Science Innovation-Driven Development of Chongqing (Grant No. 2023NSCQ-LZX0218), Sichuan Science and Technology Program (Grant No. 2023NSFSCO060).
	
	\section*{Data availability statement}
	Our code is available at \url{https://github.com/fsliuxl/LGPM-ISTA}
	
	\bibliographystyle{abbrv}
	\bibliography{l1tv}
	
\end{document}